\documentclass{article} 

\usepackage{latexsym}
\usepackage{amsfonts}
\usepackage{amssymb}
\usepackage{amsmath}
\usepackage{amscd}
\usepackage{eucal}
\usepackage[all, knot]{xy}
\xyoption{arc}
\xyoption{rotate}
\xyoption{dvips}
\usepackage{amsthm}
\usepackage[dvips]{pict2e}

\usepackage{verbatim} 
\usepackage[dvips]{graphicx}

\numberwithin{equation}{section}

\newtheorem{Def}{Definition}[section]
\newtheorem{thm}[Def]{Theorem}

\newtheorem{prop}[Def]{Proposition}
\newtheorem{cor}[Def]{Corollary}

\newtheoremstyle{example}{\topsep}{\topsep}%
     {}
     {}
     {\bfseries}
     {.}
     {8pt}
     {\thmname{#1}\thmnumber{ #2}\thmnote{ #3}}

   \theoremstyle{example}
   
\newtheorem{rem}[Def]{Remark}
\newtheorem{eg}[Def]{Example}

\newcommand{\glob}{
\xy
(-5,0)*+{.}="1";
(5,0)*+{.}="2";
{\ar@/^1pc/^{} "1";"2"};
{\ar@/_1pc/_{} "1";"2"};
{\ar@{=>}^{} (0,2)*{};(0,-2)*{}} ;
\endxy}

\newcommand{\lra}{\ensuremath{\longrightarrow}}
\newcommand{\cat}[1]{\ensuremath{\mbox{\bfseries {\upshape {#1}}}}}
\newcommand{\cl}[1]{\ensuremath{\mathcal {#1}}}

\newcommand{\ed}{\end{document}}
\newcommand{\map}[1]{\ensuremath{\stackrel{{#1}}{\lra}}}

\newcommand{\set}{\cat{Set}}

\newcommand{\Cat}{{\cat{Cat}}}

\newcommand{\numarabic}{\renewcommand{\labelenumi}{\arabic{enumi}.}}

\newcommand{\psinv}[6]{\xy
(#1,0)*+{#3}="x";
(#2,0)*+{#5}="y";
{\ar@<.7ex>^{#4} "x"; "y"};
{\ar@<.7ex>^{#6} "y"; "x"};
\endxy
}

\newcommand{\mnd}[1]{\ensuremath{\mbox{{\upshape Mnd}}({#1}) }}
\newcommand{\mndnb}[1]{\ensuremath{\mbox{{\upshape Mnd}}^{#1}(\cl{B}) }}
\newcommand{\mndb}{\ensuremath{\mbox{{\upshape Mnd}}(\cl{B}) }}

\begin{document}

\title{Iterated distributive laws} 

\author{Eugenia Cheng\\Department of Mathematics, Universit\'e de Nice Sophia-Antipolis \\ and \\ Department of Mathematics, University of Sheffield \\E-mail: e.cheng@sheffield.ac.uk}

\date{October 2007}

\maketitle

\begin{abstract}
We give a framework for combining $n$ monads on the same category via distributive laws satisfying Yang-Baxter equations, extending the classical result of Barr and Wells which combines two monads via one distributive law.  We show that this corresponds to iterating $n$-times the process of taking the 2-category of monads in a 2-category, extending the result of Street characterising distributive laws.  We show that this framework can be used to construct the free strict $n$-category monad on $n$-dimensional globular sets; we first construct for each $i$ a monad for composition along bounding $i$-cells, and then we show that the interchange laws define distributive laws between these monads, satisfying the necessary Yang-Baxter equations. 
\end{abstract}

%

\setcounter{tocdepth}{1}

\tableofcontents

\section*{Introduction}

\addcontentsline{toc}{section}{Introduction}

Monads give us a way of describing algebraic structures such as monoids, groups, rings and categories.  Distributive laws between monads give us a way of combining two such algebraic structures in a coherent way.  For example, a ring is a monoid under multiplication and an abelian group under multiplication, where the mulplication and addition must interact coherently.  We can thus construct the monad for rings using the monad for monoids and the monad for abelian groups, via a distributive law.  In this work we describe how to extend this to combine three or more algebraic structures in a coherent way.  Our motivating example is the construction of strict $n$-categories; other examples are rings, constructed using three monads instead of the usual two, and rigs (semi-rings), constructed using four.  

The starting point for distributive laws is that we have two monads $S$ and $T$, say, on a category \cl{C}.  We might then want to ask: can we compose them to get a monad $TS$?  If we think of the action of a monad as building in algebraic structure freely then $TS$ would first build in $S$-like structure freely, and then $T$-like structure ``on top''.  For this to be a monad we can ask that the $S$-like structure ``distribute'' over the $T$-like structure, via a natural transformation
\[ ST \stackrel{\lambda}{\Rightarrow} TS\]
which we can think of as ``moving the $S$ structure across the $T$-structure''.   We can then apply the multiplication for $S$ and $T$ to get a putative multiplication for $TS$ 
\[\xy
(40,0)*+{TS}="2"; 
(22,0)*+{T^2 S^2}="1"; 
(0,0)*+{TSTS}="0"; 
{\ar@{=>}^{T \lambda S} "0"; "1"};
{\ar@{=>}^>>>>>{\mu^T \mu^S} "1"; "2"};
\endxy;\]
of course, we then need to check that this satisfies the associativity axiom for a monad.  

A \emph{distributive law} of $S$ over $T$ is defined in \cite{bec1} to be a natural transformation $\lambda$ as above, satisfying axioms ensuring coherent interaction with the monad structures for $S$ and $T$.  One consequence is that the putative multiplication we defined above really does make $TS$ into a monad (with unit $\eta^T \eta^S :1 \Rightarrow TS$); note that distributive laws are directional, and $ST$ does not become a monad.

We can now examine the situation involving three monads $A$, $B$ and $C$, say, on the same category, with distributive laws
\[
\begin{array}{ccc}
\lambda &:& BA \Rightarrow AB\\
\sigma &:& CA \Rightarrow AC\\
\tau &:& CB \Rightarrow BC.
\end{array}
\]
By the above distributive law result, we have canonical monads $AB$, $AC$ and $BC$; we can further ask if we have a monad $ABC$.  This would certainly follow from the theory of distributive laws if we had a distributive law
\[ (BC)A \Rightarrow A(BC)\]
or indeed
\[C(AB) \Rightarrow (AB)C\]
each of which would give a canonical monad $ABC$.  It turns out that although we can easily construct such maps, they will not automatically be distributive laws.  However we can ensure that they are distributive laws by imposing the following axiom:
\[\xy
(0,12)*+{C B A}="1";
(12,24)*+{B C A}="2";
(36,24)*+{B A C}="3";
(12,0)*+{C A B}="4";
(36,0)*+{A C B}="5";
(48,12)*+{A B C}="6";
{\ar^{\lambda A} "1";"2"};
{\ar^{B \sigma} "2";"3"};
{\ar^{\tau C} "3";"6"};
{\ar_{C\tau} "1";"4"};
{\ar_{\sigma B} "4";"5"};
{\ar_{A \lambda} "5";"6"};
\endxy\]
If we examine how $A$, $B$, and $C$ are ``switched'' past each other in this diagram, we see that it is none other than the Yang-Baxter equation.  This turns out to be all we need to make sure we induce the two new distributive laws above.  Moreover, the two resulting monad structures on $ABC$ coincide.

More generally we can consider a series of $n$ monads $T_1, \ldots T_n$ on a category, with pairwise distributive laws going in one ``direction'' only, that is, for all $i > j$ a distributive law
\[T_i T_j \Rightarrow T_j T_i\]
where every three monads satisfies the Yang-Baxter equation; we state this precisely in the main theorem (Theorem~\ref{mainthm}) and call it a ``distributive series of monads''.  This turns out to be enough to ensure that the monads are as coherent as we want --- they can be composed (in fixed order) to produce one combined monad.  There are no further conditions that need to be imposed other than the three-way Yang-Baxter conditions.  This is the main theorem of this work.  Note that this idea is essentially the same as the fact that, to express the braid group by generators and relations, a single crossing and the Yang-Baxter equation suffice.

This theorem can be proved directly by a diagram chase and induction, and this is what we do in Section~\ref{iterate}.  In Section~\ref{aussie} we look at a more abstract approach to monads provided by Street in \cite{str1}.  Usually monads are defined by giving a category \cl{C}, a functor $T: \cl{C} \lra \cl{C}$ and natural transformations 
\[\xy
(26,0)*+{T^2}="2"; 
(13,0)*+{T}="1"; 
(0,0)*+{1}="0"; 
{\ar@{=>}^{\eta} "0"; "1"};
{\ar@{=>}_>>>>>>{\mu^2} "2"; "1"};
\endxy\]
satisfying the usual unit and associativity axioms.  Of course, we could make the exact same definition inside any 2-category --- here we happen to have done it inside \cat{Cat}, the 2-category of categories, functors and natural transformations.  Moreover, Street also defines monad functors and monad transformations, so that given any 2-category \cl{B} there is a whole 2-category $\mbox{Mnd}(\cl{B})$ of monads in \cl{B}.  So we can now iterate this and take monads in $\mbox{Mnd}(\cl{B})$.  

At the end of \cite{str1} Street shows that monads in $\mbox{Mnd}(\cl{B})$ are precisely distributive laws.  That is, a 0-cell in $\mnd{\mnd{\cl{B}}}$ is given by a 0-cell in \cl{B}, two monads on it, and a distributive law of one over the other.  In Section~\ref{iteratemnd} we take Street's construction further and construct the entire 2-category $\mndnb{2} = \mnd{\mnd{\cl{B}}}$, enabling us to iterate $n$ times and get the 2-category $\mndnb{n}$.  The main theorem is then that an object in \mndnb{n} is precisely a distributive series of $n$ monads as in the main theorem described above.  The proof is by induction and hinges on a careful (and notationally fiddly) characterisation of the 1- and 2-cells of \mndnb{n}.  

We end this section with some brief remarks on two other characterisations of monads that may shed some light on this situation: the simplicial resolution of a monad, and monads in \cl{B} via the ``free living monad'' 2-category $\Delta$.  A monad in \cl{B} is precisely a 2-functor $\Delta \lra \cl{B}$.  We then use the fact \cite{str1} that Mnd is itself a monad (on the category \cat{2-Cat} of 2-categories and 2-functors).  Thus it can be expressed as a functor $\Delta \lra \cat{Cat}$.  Using the closed structure of \cat{2-Cat} with respect to the lax Gray tensor product (see \cite{gra1, gra2}), we see that 
\[\mndnb{n} \simeq [\Delta^{\otimes n}, \cl{B}] \]
and the Yang-Baxter equations for the monads correspond to Yang-Baxter equations arising from the relations on the generators defining the Gray tensor product.  

In Section~\ref{interchange} we present our motivating example, the construction of the free strict $n$-category monad (on $n$-dimensional globular sets) by building up the monad from $n$ separate monads for composition.  That is, we isolate composition along $i$-cell boundaries for each $0 \leq i \leq n-1$, and construct for each $i$ a monad $T_i$ that gives this composition alone.  These monads are simply a generalisation of the usual free category monad on graphs.  It is then an interesting fact that the required distributive laws come from the interchange laws for a strict $n$-category -- for all $i<j$ there is a strict interchange law for $j$-composition over $i$-composition generalising the usual interchange law in 2-categories, and it turns out that this does indeed define a distributive law of the monad $T_j$ over the monad $T_i$ where $i<j$.  Moreover, these distributive laws satisfy the Yang-Baxter equation for all $i<j<k$.  Thus, using our main theorem we can construct the free strict $n$-category monad as the composite of the monads $T_i$, in order.  

This is reminiscent of the definition of strict $n$-category as an $n$-globular set in which every sub-2-globular set is a 2-category, that is, where interchange holds for every pair of dimensions (see \cite{str2}).  

It is worth commenting briefly on the notion of strict interchange.  Interchange in $n$-categories is sometimes thought of as ``the only part of weak $n$-categories that cannot be strictified''.  More precisely, we know that not every tricategory is triequivalent to a completely strict 3-category; the well-known coherence result is that every tricategory is triequivalent to a Gray-category \cite{gps1, gur2}.  This has led to a generally accepted conjecture that every weak $n$-category should be equivalent to a semi-strict one, where ``semi-strict'' means that every thing is strict \emph{apart} from interchange -- interchange thus being the only part that cannot be strictified.  However, recent work of Joyal and Kock \cite{jk1} has shown that a different type of ``semi-strict 3-category'' is also fruitful: one in which everything is strict apart from units.  Joyal and Kock have shown that such 3-categories do model homotopy 3-types, and do give rise to braided monoidal categories in the suitably degenerate case, two crucial things that fully strict 3-categories cannot do.  Indeed Simpson \cite{sim1, sim2} conjectures that $n$-categories with weak units (and everything else strict) are enough to model $n$-types for all $n$.  

In the light of these recent results we consider the study of $n$-categories with strict interchange to be important, and the results of the present work will contribute towards that study.  In a future work, and in \cite{che16} we will use the theory of iterated distributive laws to study Trimble's notion of $n$-category \cite{tri1}, which has strict interchange but weak units.  The results of the present work enable us to construct a comparison functor from Trimble $n$-categories to those defined by Batanin \cite{bat1}.    

Finally we note that the proofs and constructions in this work are all completely straightforward (in that there are no surprises) though sometimes lengthy and afflicted with unavoidably complex notation.

\subsubsection*{Acknowledgements}

I would like to thank Michael Shulman for useful conversations relating to this work.

\section{Distributive laws}\label{classical}

We first recall the classical theory of distributive laws.  

\begin{Def} (Beck \cite{bec1})

Let $S$ and $T$ be monads on a category \cl{C}.  A {\em distributive law of $S$ over $T$} consists of a natural transformation $\lambda : ST \Rightarrow TS$ such that the following diagrams commute.

\begin{equation}\label{dl1}
\xymatrix{
& T \ar[dl]_{\eta^S T} \ar[dr]^{T\eta^S} && S^2 T \ar[r]^<<<<{S\lambda} \ar[d]_{\mu^S T}
& STS \ar[r]^{\lambda S }
& TS^2 \ar[d]^{T \mu^S} \\
ST \ar[rr]_{\lambda} && TS & ST \ar[rr]_{\lambda}
&& TS}
\end{equation}

\begin{equation}\label{dl2}
\xymatrix{
& S \ar[dl]_{S\eta^T} \ar[dr]^{\eta^T S} && ST^2 \ar[r]^<<<<{\lambda T} \ar[d]_{S\mu^T}
& TST \ar[r]^{T\lambda}
& T^2 S \ar[d]^{\mu^T S} \\ 
ST \ar[rr]_{\lambda} && TS & ST \ar[rr]_{\lambda}
&& TS }
\end{equation}

\end{Def}

\noindent Note that the first pair of axioms (\ref{dl1}) is telling us that ``$\lambda$ interacts coherently with the monad structure for $S$'' and the second pair (\ref{dl2}) that ``$\lambda$ interacts coherently with the monad structure for $T$''.  

The main theorem about distributive laws tells us about new monads that arise canonically as a result of the distributive law.  

\begin{thm}[Barr and Wells, \cite{bw1}]\label{barrwells}

The following are equivalent:
\numarabic
\begin{itemize}
\item A distributive law of $S$ over $T$.
\item A lifting of the monad $T$ to a monad $T'$ on $S\mbox{{\upshape -Alg}}$.
\item A lifting of the monad $S$ to a monad $S'$ on $\mbox{{\upshape Kl}} (T)$.
\end{itemize}
It follows that $TS$ canonically acquires the structure of a monad, whose category of algebras coincides with that of the lifted monad $T'$, and whose Kleisli category coincides with that of $S'$.
\end{thm}

\subsubsection*{Remark}

We will not be using the Kleisli part of this theorem in this work.

\bigskip

The idea is that $T$ lifts to a monad on $S$-algebras via $\lambda$, with its action on an $S$-algebra $\theta$ given by
\[\xy
{\ar^{\theta} (0,20)*+{SA};(0,10)*+{A}};
{\ar@{|->} (7,15);(14,15)};
{\ar^{\lambda_A} (23,25)*+{STA};(23,15)*+{TSA}};
{\ar^{T\theta} (23,15)*+{TSA};(23,5)*+{TA}};
\endxy;\]

\noindent $TS$ becomes a monad with multiplication 

\[\xy
(40,0)*+{TS}="2"; 
(22,0)*+{T^2 S^2}="1"; 
(0,0)*+{TSTS}="0"; 
{\ar@{=>}^{T \lambda S} "0"; "1"};
{\ar@{=>}^>>>>>{\mu^T \mu^S} "1"; "2"};
\endxy;\]
and unit $\eta^T \eta^S : 1 \Rightarrow TS$. The axioms for $\lambda$ precisely ensure that both of these satisfy the axioms for a monad.

\begin{eg}\label{ring}

{\bfseries (Rings)}

$\cl{C} = \cat{Set}$

$S = \mbox{free commutative monoid monad}$

$T = \mbox{free abelian group monad}$

$\lambda = \mbox{the usual distributive law for multiplication and addition}$ eg \[(a+b)(c+d) \mapsto ac+bc+ad+bd.\]  

Then the composite monad $TS$ is the free ring monad.
\end{eg}

\begin{eg}\label{monoidpointing}

{\bfseries (Monoids)}

$\cl{C}=\cat{Set}$

$S = \mbox{monad for non-unital associative multiplication}$

$T = \mbox{monad for pointed sets}$ ie $TA = A \coprod \{\ast\}$

$\lambda$ ensures that $\ast$ acts as a unit for multiplication:
\[\begin{array}{ccc}
S(A \coprod \{\ast\}) & \lra & SA \coprod \{\ast\} \\
a_1 . \ldots . a_i . \ast . a_{i+1} . \ldots . a_n & \mapsto & a_1 . \ldots . a_n \end{array}\]

Then the composite monad $TS$ is the free monoid monad.
\end{eg}  

\begin{eg}\label{commonoidpointing}

{\bfseries (Commutative monoids)}

As above but with $S$ giving a commutative multiplication; note that this does not work for groups.
\end{eg}

\begin{eg}

{\bfseries (2-categories)}

$\cl{C}=\cat{2-GSet}$, the category of 2-globular sets, that is diagrams in \cat{Set}
\[\xy
(120,0)*+{A(0)}="0"; 
(100,0)*+{A(1)}="1"; 
(80,0)*+{A(2)}="2"; 
{\ar@<0.7ex>^{s} "1"; "0"};
{\ar@<-0.7ex>_{t} "1"; "0"};
{\ar@<0.7ex>^{s} "2"; "1"};
{\ar@<-0.7ex>_{t} "2"; "1"};
\endxy\]

\noindent with $ss=ts, ts=tt$.

$S = $ monad for vertical composition of 2-cells ($A(1)$ and $A(0)$ are unchanged)

$T=$ monad for horizontal composition of 2-cells and 1-cells ($A(0)$ is unchanged)

$\lambda$ is given by the interchange law eg
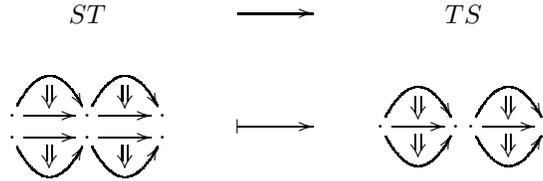
\[\xy
(0,15)*{ST};
(50,15)*{TS};
{\ar (20,15);(30,15)};
{\ar@{|->} (20,0);(30,0)};
(50,0)*{ 	
\xy
(0,0)*{
\xy
(0,0)*+{.}="1";
(10,0)*+{.}="2";
{\ar@/^1.25pc/^{} "1";"2"};
{\ar@/_1.25pc/_{} "1";"2"};
{\ar "1";"2"};
{\ar@{=>}^{} (5,4)*{};(5,1)*{}} ;
{\ar@{=>}^{} (5,-1)*{};(5,-4)*{}} ;
\endxy};
(12,0)*{
\xy
(0,0)*+{.}="1";
(10,0)*+{.}="2";
{\ar@/^1.25pc/^{} "1";"2"};
{\ar@/_1.25pc/_{} "1";"2"};
{\ar "1";"2"};
{\ar@{=>}^{} (5,4)*{};(5,1)*{}} ;
{\ar@{=>}^{} (5,-1)*{};(5,-4)*{}} ;
\endxy}
\endxy};
(0,0)*{ 
\xy
(0,0)*{
\xy
(0,0)*+{.}="1";
(10,0)*+{.}="2";
(20,0)*+{.}="3";
{\ar@/^1.25pc/^{} "1";"2"};
{\ar "1";"2"};
{\ar@{=>}^{} (5,4)*{};(5,1)*{}} ;
{\ar@/^1.25pc/^{} "2";"3"};
{\ar "2";"3"};
{\ar@{=>}^{} (15,4)*{};(15,1)*{}} ;
\endxy};
(0,-8)*{ 
\xy
(0,0)*+{.}="1";
(10,0)*+{.}="2";
(20,0)*+{.}="3";
{\ar@/_1.25pc/^{} "1";"2"};
{\ar "1";"2"};
{\ar@{=>}^{} (5,-1)*{};(5,-4)*{}} ;
{\ar@/_1.25pc/^{} "2";"3"};
{\ar "2";"3"};
{\ar@{=>}^{} (15,-1)*{};(15,-4)*{}} ;
\endxy};
\endxy};
\endxy\]

\noindent Note that this distributive law can only go in this direction; we will discuss this in more detail in Section~\ref{interchange}.  The fact that this is a distributive law can be proved by direct calculation; alternatively it follows from abstract results that we give in Theorem~\ref{fourpointseven}.

\end{eg}

\section{Iterated distributive laws}\label{iterate}

In this section we generalise the notion of distributive law to the case when we have more than two monads interacting with each other.

\subsection{The main theorem}

\begin{thm}\label{mainthm}

Fix $n \geq 3$.  Let $T_1, \ldots, T_n$ be monads on a category \cl{C}, equipped with

\begin{itemize}

\item for all $i > j$ a distributive law $\lambda_{ij}: T_i T_j \Rightarrow T_j T_i$, satisfying

\item for all $i>j>k$ the ``Yang-Baxter'' equation given by the commutativity of the following diagram

\begin{equation}\label{yb}
\xy
(0,12)*+{T_i T_j T_k}="1";
(12,24)*+{T_j T_i T_k}="2";
(36,24)*+{T_j T_k T_i}="3";
(12,0)*+{T_i T_k T_j}="4";
(36,0)*+{T_k T_i T_j}="5";
(48,12)*+{T_k T_j T_i}="6";
{\ar^{\lambda_{ij} T_k} "1";"2"};
{\ar^{T_j \lambda_{ik}} "2";"3"};
{\ar^{\lambda_{jk} T_i} "3";"6"};
{\ar_{T_i\lambda_{jk}} "1";"4"};
{\ar_{\lambda_{ik} T_j} "4";"5"};
{\ar_{T_k \lambda_{ij}} "5";"6"};
\endxy\end{equation}

%
%

\end{itemize}

\noindent Then for all $1\leq i < n$  we have canonical monads 
\[T_1 T_2 \cdots T_i \quad \mbox{and} \quad T_{i+1} T_{i+2} \cdots T_n\]

\noindent together with a distributive law of \ $T_{i+1} T_{i+2} \cdots T_n$\  over \ $T_1 T_2 \cdots T_i$\ i.e. 
\[(T_{i+1} T_{i+2} \cdots T_n) (T_1 T_2 \cdots T_i) \Rightarrow (T_1 T_2 \cdots T_i)(T_{i+1} T_{i+2} \cdots T_n)\]

\noindent given by the obvious composites of the $\lambda_{ij}$.  Moreover, all the induced monad structures on \ $T_1 T_2 \cdots T_n$ are the same.

\end{thm}

\begin{rem}\label{mainthmrem}

In some situations it may be convenient to index the monads in the opposite direction.  For example, the direction given above is more convenient for the constructions in Section~\ref{interchange}, but the opposite direction is needed in the proofs in Section~\ref{aussie}.  This does not of course affect the content of the theorem, just the notation.  Reversing the indexing gives a series $T_1, \ldots, T_n$ of monads on a category \cl{C}, equipped with

\begin{itemize}

\item for all $i < j$ a distributive law $\lambda_{ij}: T_i T_j \Rightarrow T_j T_i$, satisfying

\item for all $i<j<k$ the ``Yang-Baxter'' diagram

\[\xy
(0,12)*+{T_i T_j T_k}="1";
(12,24)*+{T_j T_i T_k}="2";
(36,24)*+{T_j T_k T_i}="3";
(12,0)*+{T_i T_k T_j}="4";
(36,0)*+{T_k T_i T_j}="5";
(48,12)*+{T_k T_j T_i}="6";
{\ar^{\lambda_{ij} T_k} "1";"2"};
{\ar^{T_j \lambda_{ik}} "2";"3"};
{\ar^{\lambda_{jk} T_i} "3";"6"};
{\ar_{T_i\lambda_{jk}} "1";"4"};
{\ar_{\lambda_{ik} T_j} "4";"5"};
{\ar_{T_k \lambda_{ij}} "5";"6"};
\endxy\]

%
%

\end{itemize}

\noindent Then for all $1\leq i < n$  we have canonical monads 
\[T_n T_{n-1} \cdots T_{i+1} \quad \mbox{and} \quad T_i T_{i-1} \cdots T_1\]

\noindent together with a distributive law of \ $T_i T_{i-1} \cdots T_1$\  over \ $T_n T_{n-1} \cdots T_{i+1}$\ i.e. 
\[(T_i T_{i-1} \cdots T_1) (T_n T_{n-1} \cdots T_{i+1}) \Rightarrow (T_n T_{n-1} \cdots T_{i+1}) (T_i T_{i-1} \cdots T_1)\]

\noindent given by the obvious composites of the $\lambda_{ij}$.  In this case, there are induced monad structures on \ $T_n T_{n-1} \cdots T_1$, and they are all the same.

\end{rem}

\begin{proof} By induction (using the notation of the Theorem, not the remark).  The case $n=3$ is a simple diagram chase as follows.  For ease of notation we write the monads $T_1, T_2, T_3$ as $A,B,C$ with distributive laws
\[
\begin{array}{ccc}
\lambda &:& BA \Rightarrow AB\\
\sigma &:& CA \Rightarrow AC\\
\tau &:& CB \Rightarrow BC.
\end{array}
\]
and we will check that the following is a distributive law:
\[\xy
(36,0)*+{ABC}="2"; 
(18,0)*+{BAC}="1"; 
(0,0)*+{BCA}="0"; 
{\ar@{=>}^{B\sigma} "0"; "1"};
{\ar@{=>}^{\tau C} "1"; "2"};
\endxy.\]
Now, the triangular axioms for a distributive law follow immediately without any need for the Yang-Baxter condition, as does the pentagon axiom (\ref{dl2}): 
%
%
\[
\xy
(27,17)*+{\mbox{{\bfseries (1.2)}}}="10";
(93,17)*+{\mbox{{\bfseries (1.2)}}}="11";
(0,40)*+{BCAA}="1";
(30,40)*+{BACA}="2";
(60,40)*+{ABCA}="3";
(90,40)*+{ABAC}="4";
(120,40)*+{AABC}="5";
(60,20)*+{BAAC}="6";
(0,0)*+{BCA}="7";
(60,0)*+{BAC}="8";
(120,0)*+{ABC}="9";
{\ar^{B\sigma A} "1";"2"};
{\ar^{\lambda CA} "2";"3"};
{\ar^{AB\sigma} "3";"4"};
{\ar^{A\lambda C} "4";"5"};
{\ar_{BA\sigma} "2";"6"};
{\ar_{\lambda AC} "6";"4"};
{\ar^{BC\mu^A} "1";"7"};
{\ar^{\mu^A BC} "5";"9"};
{\ar^{B\sigma} "7";"8"};
{\ar^{\lambda C} "8";"9"};
{\ar^{B\mu^AC} "6";"8"};
\endxy
\]
where the unmarked region commutes by naturality.  The pentagon axiom (\ref{dl1}) is seen to commute as follows:

\vspace{5ex}


\hspace*{-6em}$\xy
(68,50)*+{\mbox{{\textsf{\bfseries Yang-Baxter}}}}="20";
(34,10)*+{{\mbox{{\bfseries (1.1)}}}}="21";
(102,10)*+{{\mbox{{\bfseries (1.1)}}}}="22";
(0,60)*+{BCBCA}="1";
(34,60)*+{BCBAC}="2";
(68,60)*+{BCABC}="3";
(102,60)*+{BACBC}="4";
(136,60)*+{ABCBC}="5";
(0,40)*+{BBCCA}="6";
(34,40)*+{BBCAC}="7";
(68,40)*+{BBACC}="8";
(102,40)*+{BABCC}="9";
(136,40)*+{ABBCC}="10";
(12,20)*+{BCCA}="11";
(34,20)*+{BCAC}="12";
(56,20)*+{BACC}="13";
(80,20)*+{BBAC}="14";
(102,20)*+{BABC}="15";
(124,20)*+{ABBC}="16";
(0,0)*+{BCA}="17";
(68,0)*+{BAC}="18";
(136,0)*+{ABC}="19";
{\ar^{BCB\sigma} "1";"2"};
{\ar^{BC\lambda C} "2";"3"};
{\ar^{B\sigma BC} "3";"4"};
{\ar^{\lambda CBC} "4";"5"};
{\ar_{B\tau CA} "1";"6"};
{\ar^{B\tau AC} "2";"7"};
{\ar^{BA\tau C} "4";"9"};
{\ar^{AB\tau C} "5";"10"};
{\ar^{BBC\sigma} "6";"7"};
{\ar^{BB\sigma C} "7";"8"};
{\ar^{B\lambda CC} "8";"9"};
{\ar^{AB\tau C} "9";"10"};
{\ar_{\mu^B \mu^C A} "6";"17"};
{\ar^{\mu^B CCA} "6";"11"};
{\ar^{\mu^B CAC} "7";"12"};
{\ar_{\mu^B ACC} "8";"13"};
{\ar^{BBA\mu^C} "8";"14"};
{\ar^{BAB\mu^C} "9";"15"};
{\ar_{ABB\mu^C} "10";"16"};
{\ar^{A\mu^B \mu^C} "10";"19"};
{\ar^{BC\sigma} "11";"12"};
{\ar^{B\sigma C} "12";"13"};
{\ar^{B\lambda C} "14";"15"};
{\ar^{\lambda BC} "15";"16"};
{\ar^{B\mu^C A} "11";"17"};
{\ar_{BA \mu^C} "13";"18"};
{\ar^{\mu^BAC} "14";"18"};
{\ar_{A\mu^B C} "16";"19"};
{\ar^{B\sigma} "17";"18"};
{\ar^{\lambda C} "18";"19"};
\endxy$

\vspace{5ex}

\noindent where all the unmarked parts commute by naturality, and the Yang-Baxter equation marked is that for $C,B,A$ with an extra $B$ applied on the left throughout, and a $C$ on the right.  The result for the other distributive law follows similarly.

Now consider $n>3$.  Given $i$ with $1 \leq i <n$, the series of monads $T_1, \ldots, T_i$ and the series $T_{i+1}, \ldots, T_n$ is each a series of monads satisfying the induction hypothesis; each series has fewer than $n$ monads so by induction we have the monads \ $T_1 T_2 \cdots T_i$\ and \ $T_{i+1} T_{i+2} \cdots T_n$ as required.  Now define monads $S_1, \ldots, S_{n-1}$ by

\begin{itemize}
\item $S_i = T_i$ for all $1 \leq i < n-1$, and
\item $S_{n-1} = T_{n-1} T_n $.  
\end{itemize}

\noindent We now check that the monads $S_1, \ldots, S_{n-1}$ satisfy the hypothesis of the theorem:

\numarabic

\begin{enumerate}
\item We need for all $i>j$ a distributive law of $S_i$ over $S_j$.  For the cases $i<n-1$ this is just the distributive law of $T_i$ over $T_j$.  For the case $i=n-1$ we need for all $j<n-1$ a distributive law of $T_{n-1} T_n$ over $T_i$; this follows from the result for $n=3$ applied to the monads $T_i, T_{n-1}, T_n$.  

\item We need for all $i>j>k$ the Yang-Baxter equation for $S_i, S_j, S_k$.  Again, for the cases $i<n-1$ this is just the Yang-Baxter eqation for $T_i, T_j, T_k$.  For the case $i=n-1$ we need for all $k<j<n-1$ the Yang-Baxter equation for the monads $S_{n-1}, T_j, T_k$, that is, the monads $(T_{n-1}T_n), T_j, T_k$.  This follows from the Yang-Baxter equations for $T_{n-1}, T_j, T_k$ and $T_n, T_j, T_k$.  This is seen by the following diagram, where for clarity we have labelled the four monads in question $A,B,C,D$:
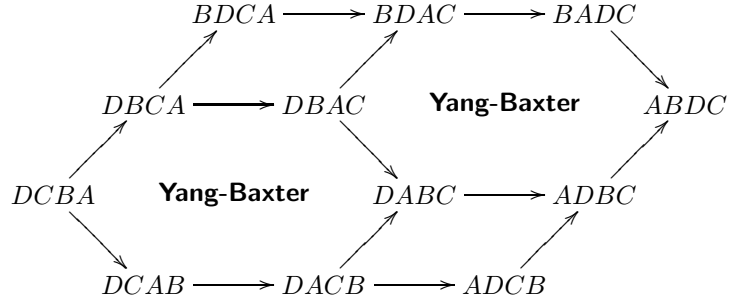
\[\xy
(24,12)*+{\mbox{{\textsf{\bfseries Yang-Baxter}}}};
(60,24)*+{\mbox{{\textsf{\bfseries Yang-Baxter}}}};
(0,12)*+{DCBA}="1";
(12,24)*+{DBCA}="2";
(36,24)*+{DBAC}="3";
(12,0)*+{DCAB}="4";
(36,0)*+{DACB}="5";
(48,12)*+{DABC}="6";
(24,36)*+{BDCA}="7";
(48,36)*+{BDAC}="8";
(72,36)*+{BADC}="9";
(84,24)*+{ABDC}="10";
(72,12)*+{ADBC}="11";
(60,0)*+{ADCB}="12";
{\ar^{} "1";"2"};
{\ar^{} "2";"3"};
{\ar^{} "3";"6"};
{\ar_{} "1";"4"};
{\ar_{} "4";"5"};
{\ar_{} "5";"6"};
{\ar_{} "2";"7"};
{\ar_{} "7";"8"};
{\ar_{} "8";"9"};
{\ar_{} "9";"10"};
{\ar_{} "3";"8"};
{\ar_{} "6";"11"};
{\ar_{} "11";"10"};
{\ar_{} "12";"11"};
{\ar_{} "5";"12"};
\endxy\]
Here the lower left hexagon is the Yang-Baxter equation for $B,C,A$ (with $D$ applied on the left), and the upper right hexagon is the Yang-Baxter equation for $D,B,A$ (with $C$ applied on the right).  

\end{enumerate}

\noindent So by the result for $n-1$ we have distributive laws for all $1\leq i < n-1$
\[(S_{i+1} S_{i+2} \cdots S_{n-1}) (S_1 S_2 \cdots S_i) \Rightarrow (S_1 S_2 \cdots S_i)(S_{i+1} S_{i+2} \cdots S_{n-1})\]

\noindent that is, 
\[(T_{i+1} T_{i+2} \cdots T_{n-1} T_n) (T_1 T_2 \cdots T_i) \Rightarrow (T_1 T_2 \cdots T_i)(T_{i+1} T_{i+2} \cdots T_{n-1} T_n)\]

\noindent all inducing the same monad. We are still missing the distributive law for $i=n-1$; for this we just repeat the above proof using monads $(T_1 T_2), T_3, \ldots, T_n$, and the full result follows.  
\end{proof}

\begin{Def}
A {\em distributive series of $n$ monads} is a system of monads and distributive laws as in Theorem~\ref{mainthm}.  
\end{Def}

\subsection{Examples}

In this section we give two brief examples of iterated distributive laws; our main example, that of interchange for $n$-categories will be treated in Section~\ref{interchange}.  


\subsubsection{Rings}

The usual example for distributive laws involves constructing the free ring monad from the free commutative monoid monad and the free abelian group monad (Example~\ref{ring}).  However, we can also construct it from the following distributive series of 3 monads on \set, combining Example~\ref{ring} with Example~\ref{commonoidpointing}:
\[\begin{array}{ccl}
A &=& \mbox{monad for associative non-unital binary multiplication $\times$} \\
B &=& \mbox{monad for pointed sets i.e. $X \mapsto X \coprod \{1\}$} \\
C &=& \mbox{free additive abelian group monad}
\end{array}\] 

\noindent We have distributive laws:

\begin{itemize}
\item $AB \Rightarrow BA$ as in Example~\ref{commonoidpointing}, ensuring that 1 acts as a unit for $\times$
\item $AC \Rightarrow CA$ as in Example~\ref{ring}, the usual distributive law of $\times$ over $+$, but without units, and
\item $BC \Rightarrow CB$ is the obvious embedding
\end{itemize}
and it is easy to check the Yang-Baxter equation.  The composite monad $CBA$ is then the free ring monad.  

\subsubsection{Rigs}

If we have no negatives we can further decompose the situation with the following four monads:
\[\begin{array}{ccl}
A &=& \mbox{monad for associative non-unital binary multiplication $\times$} \\
B &=& \mbox{monad for pointed sets i.e. $X \mapsto X \coprod \{1\}$} \\
C &=& \mbox{monad for associative non-unital non-invertible $+$}\\
D &=& \mbox{$X \mapsto X \coprod \{0\}$}
\end{array}\] 
and we have distributive laws as follows:

\begin{itemize}
\item $AB \Rightarrow BA$ as above,
\item $AC \Rightarrow CA$ as above but without 0,
\item $AD \Rightarrow DA$ ensures that $x \times 0 = 0$,
\item $BC \Rightarrow CB$ as above but without 0,
\item $BD \Rightarrow DB$ is given by the obvious map $X \coprod \{0\} \coprod \{1\} \lra X \coprod \{1\} \coprod \{0\}$
\item $CD \Rightarrow DC$ ensures that $x+0 = 0$
\end{itemize}
Then we can check all the required Yang-Baxter equations, and the resulting composite monad $DCBA$ is the free rig monad.

\section{Iterated distributive laws via the formal theory of monads}\label{aussie}

In his classic paper {\em The formal theory of monads} \cite{str1} Street defines for any 2-category \cl{B} a 2-category $\mbox{Mnd}(\cl{B})$ of monads in \cl{B}.  Then distributive laws arise as monads in $\mbox{Mnd}(\cl{B})$.  In this section we show that iterating this process gives the distributive series of monads described in the previous section.

\subsection{The 2-category of monads in \cl{B}}\label{ftm}

First we recall the basic definitions as given in \cite{str1}; none of the material in this section is new.

\begin{Def}
Let \cl{B} be a 2-category.  A monad in \cl{B} is given by

\begin{itemize}

\item a 0-cell $X$, 
\item a 1-cell $S:X \lra X$, and 
\item 2-cells $1 \stackrel{\eta}{\Rightarrow} S \stackrel{\mu}{\Leftarrow} S^2$, satisfying

\item the usual monad axioms

\begin{equation}\label{monadax}
\xymatrix{
S^2 \ar[r]^{S\eta} \ar[dr]_{1} & S^2 \ar[d]^{\mu} & S \ar[l]_{\eta S} \ar[dl]^{1} && S^3 \ar[r]^{S\mu} \ar[d]_{\mu S} & S^2 \ar[d]^{\mu}\\
& S &&& S^2 \ar[r]_{\mu} & S}\end{equation}

\end{itemize}


\noindent A monad morphism $(X, S) \lra (X', S')$ is given by 

\begin{itemize}

\item a 1-cell $T : X \lra X'$ and
\item a 2-cell $\lambda: S' T \Rightarrow T S$ i.e. 
\[\xy
(0,0)*+{X}="1";
(0,20)*+{X}="2";
(20,20)*+{X'}="3";
(20,0)*+{X'}="4";
{\ar^{T} "2";"3"};
{\ar_{T} "1";"4"};
{\ar_{S} "2";"1"};
{\ar^{S'} "3";"4"};
{\ar@{=>}^{\lambda} (12,12);(8,8)};
\endxy\]

\noindent such that

\item the following diagrams commute:

\begin{equation}\label{mmorph}
\xymatrix{
& T \ar[dl]_{\eta^{S'} T} \ar[dr]^{T\eta^{S}} && {{S'}^2} T \ar[r]^<<<<{S'\lambda} \ar[d]_{\mu^{S'} T}
& S'TS \ar[r]^{\lambda S }
& T{S}^2 \ar[d]^{T \mu^{S}} \\
S'T \ar[rr]_{\lambda} && TS & S'T \ar[rr]_{\lambda}
&& TS }
\end{equation}

\end{itemize} 

\noindent A monad transformation \[
\xy
(-12,0)*+{(X, S)}="4";
(12,0)*+{(X', S')}="6";
{\ar@/^1.65pc/^{(T, \lambda)} "4";"6"};
{\ar@/_1.65pc/_{(T', \lambda')} "4";"6"};
{\ar@{=>}^<<<{\scriptstyle \alpha} (0,3)*{};(0,-3)*{}} ;
\endxy
\]

\noindent is given by 

\begin{itemize}
\item a 2-cell $\alpha : T \Rightarrow T'$, such that
\item the following diagram commutes

\begin{equation}\label{mtrans}
\xy
(0,0)*+{S' T'}="1";
(0,20)*+{S' T}="2";
(20,20)*+{T S}="3";
(20,0)*+{T' S}="4";
{\ar^{\lambda} "2";"3"};
{\ar_{\lambda'} "1";"4"};
{\ar_{S' \alpha} "2";"1"};
{\ar^{\alpha S} "3";"4"};
\endxy
\end{equation}

\end{itemize}

\noindent Furthermore, these data organise themselves into a 2-category $\mbox{{\upshape Mnd}}(\cl{B})$ as follows:

\begin{itemize}
\item 1-cell composition 
\[\xy
(0,0)*+{(X, S)}="1";
(27,0)*+{(X', S')}="2";
(54,0)*+{(X'', S'')}="3";
(80,0)*+{(X, S)}="4";
(120,0)*+{(X'', S'')}="5";
(67,0)*{=};
{\ar^{(T, \lambda)} "1";"2"};
{\ar^{(T', \lambda')} "2";"3"};
{\ar^{(T'T,\ T' \lambda \circ \lambda' T)} "4";"5"};
\endxy\]

\item 1-cell identity $(1,1)$

\item 2-cell composition is inherited from \cl{B}.
\end{itemize}

\end{Def}

\noindent We then of course have the notion of distributive law between monads in any 2-category, the classical distributive laws being those in the 2-category \cat{Cat} of categories, functors and natural transformations.  The following theorem of Street characterises distributive laws abstractly.

\begin{thm}[Street \cite{str1}]\label{dist1}
A monad in $\mbox{{\upshape Mnd}}(\cl{B})$ is a 0-cell $X$, a pair of monads $S$ and $T$ and a distributive law $\lambda: ST \Rightarrow TS$.
\end{thm}

\begin{proof}{\em (Sketch)}  A monad in \mnd{\cl{B}} consists of

\begin{itemize}
\item a 0-cell $(X,S)$, 
\item a 1-cell $(X,S) \map{(T,\lambda)} (X,S)$, and

\item 2-cells $\xy
(52,0)*+{(T,\lambda)^2 = (T^2,\ T\lambda \circ \lambda T)}="2"; 
(18,0)*+{(T,\lambda)}="1"; 
(0,0)*+{(1,1)}="0"; 
{\ar@{=>}^{\eta} "0"; "1"};
{\ar@{=>}_>>>>>>{\mu} "2"; "1"};
\endxy$

\end{itemize}
satisfying axioms.  Hence {\em a priori} we have a monad $S$, an endomorphism $T$, and a 2-cell $\lambda: ST \Rightarrow TS$.  Now $\eta$ and $\mu$ are 2-cells of \mnd{\cl{B}} so are given by 2-cells of \cl{B}; these together with the axioms making $(T,\lambda)$ a monad in \mnd{\cl{B}} make $T$ into a monad in \cl{B}.  The axioms for $(T,\lambda)$ being a 1-cell of \mnd{\cl{B}} give the first two axioms (\ref{dl1}) for a distributive law (interaction with the monad structure of $S$) and the axioms for $\eta$ and $\mu$ to be 2-cells of \mnd{\cl{B}} give the second two axioms (\ref{dl2}) for a distributive law (interaction with the monad structure of $T$). \end{proof}

\subsection{Iterating the $\mbox{\textnormal{Mnd}}({\cl{B}})$ construction.}\label{iteratemnd}

We now show how to iterate the \mnd{\cl{B}} construction.  We will write \mndnb{2} for \mnd{\mnd{\cl{B}}}.  Theorem~\ref{dist1} tells us what the 0-cells of \mndnb{2} are; we now characterise the rest of this 2-category.  We will spell out the details as we will need them later when we characterise \mndnb{n}.  This proof is nothing more than a careful application of the definitions.

\begin{prop}
A 1-cell 
\[((X, S), (T, \lambda)) \lra ((X', S'), (T', \lambda')) \in \mndnb{2}\] 
consists of:

\begin{itemize}

\item a 1-cell $U:X \lra X' \in \cl{B}$, and 
\item  2-cells $\sigma : S' U \Rightarrow U S$ and $\tau : T' U \Rightarrow U T \in \cl{B}$, such that

\item making the following diagrams commute:

\begin{equation}\label{mmmorph1}
\xymatrix{
& U \ar[dl]_{\eta^{S'} U} \ar[dr]^{U\eta^{S}} && {S'}^2 U \ar[r]^<<<<{S'\lambda} \ar[d]_{\mu^{S'} U}
& S'US \ar[r]^{\lambda S }
& U{S}^2 \ar[d]^{U \mu^{S}} \\
S'U \ar[rr]_{\lambda} && US & S'U \ar[rr]_{\lambda}
&& US }\end{equation}

\begin{equation}\label{mmmorph2}
\xymatrix{
& U \ar[dl]_{\eta^{T'} U} \ar[dr]^{U\eta^{T}} && {T'}^2 U \ar[r]^<<<<{T'\lambda} \ar[d]_{\mu^{T'} U}
& T'UT \ar[r]^{\lambda T }
& U{T}^2 \ar[d]^{U \mu^{T}} \\
T'U \ar[rr]_{\lambda} && UT & T'U \ar[rr]_{\lambda}
&& UT }\end{equation}

\begin{equation}\label{mmmorph3}
\xy
(0,12)*+{S' T' U}="1";
(12,24)*+{T' S' U}="2";
(36,24)*+{T' U S}="3";
(12,0)*+{S' U T}="4";
(36,0)*+{U S T}="5";
(48,12)*+{U T S}="6";
{\ar^{\lambda' U} "1";"2"};
{\ar^{T' \sigma} "2";"3"};
{\ar^{\tau S} "3";"6"};
{\ar_{S' \tau} "1";"4"};
{\ar_{\sigma T} "4";"5"};
{\ar_{U \lambda} "5";"6"};
\endxy\end{equation}

\end{itemize}

\end{prop}

\noindent Note that diagrams (\ref{mmmorph1}) make $(U,\sigma)$ into a monad map $(X, S) \lra (X', S')$, and diagrams (\ref{mmmorph2}) make $(U,\tau)$ into a monad map $(X, T) \lra (X', T')$.  Diagram (\ref{mmmorph3}) is going to become the Yang-Baxter equation.

\begin{proof}
{\em A priori} a 1-cell as required consists of 
\begin{itemize}
\item a 1-cell $(U,\sigma): (X, S) \lra (X', S') \in \mnd{\cl{B}}$, and
\item a 2-cell 
\[\xy
(0,0)*+{(X,S)}="1";
(0,25)*+{(X,S)}="2";
(30,25)*+{(X',S')}="3";
(30,0)*+{(X',S')}="4";
{\ar^{(U,\sigma)} "2";"3"};
{\ar_{(U,\sigma)} "1";"4"};
{\ar_{(T,\lambda)} "2";"1"};
{\ar^{(T',\lambda')} "3";"4"};
{\ar@{=>}^{\tau} (18,15);(13,10)};
\endxy\]
\end{itemize}

\noindent such that the following diagrams commute in \mnd{\cl{B}}:
\[\xymatrix{
& (U,\sigma) \ar[dl]_{\eta.(U,\sigma)} \ar[dr]^{(U,\sigma).\eta}\\
(T', \lambda') \circ (U,\sigma) \ar[rr]_{\tau} && (U,\sigma) \circ (T, \lambda) }\]

\begin{equation}\label{diag1}
\xymatrix{
(T', \lambda')^2 \circ (U,\sigma) \ar[rr]^<<<<<<<<<<{(T', \lambda').\tau} \ar[dd]_{\mu.(U,\sigma)}
&& (T',\lambda') \circ (U,\sigma) \circ (T, \lambda) \ar[rr]^>>>>>>>>>>{\tau.(T,\lambda)}
&& (U,\sigma) \circ (T, \lambda)^2 \ar[dd]^{(U,\sigma).\mu} \\ \\
(T', \lambda') \circ (U,\sigma) \ar[rrrr]_{\tau}
&&&& (U,\sigma) \circ (T, \lambda) }\end{equation}

\noindent Now, from the 1-cell $(U,\sigma)$ we get the required cells $U$ and $\sigma$ satisfying diagrams (\ref{mmmorph1}).  Calculating the composites giving the boundaring of $\tau$ we find that $\tau$ has the form
\[
\xy
(-16,0)*+{(X, S)}="4";
(16,0)*+{(X', S')}="6";
{\ar@/^1.65pc/^{(T' U,\ T' \sigma \circ \lambda' U)} "4";"6"};
{\ar@/_1.65pc/_{(UT,\ U\lambda \circ \sigma T)} "4";"6"};
{\ar@{=>}^<<<{\scriptstyle \tau} (0,3)*{};(0,-3)*{}} ;
\endxy
\]

\noindent thus we get the required 2-cell $\tau$ satisfying diagram (\ref{mmmorph3}).  Finally diagrams (\ref{diag1}) in \mnd{\cl{B}} give us diagrams (\ref{mmmorph2}) in \cl{B} as required.  \end{proof}

\begin{prop}
A 2-cell in \mndnb{2} 	
\[
\xy
(-20,0)*+{((X, S), (T, \lambda))}="4";
(20,0)*+{((X', S'), (T', \lambda'))}="6";
{\ar@/^1.9pc/^{((U, \sigma), \tau)} "4";"6"};
{\ar@/_1.9pc/_{((U', \sigma'), \tau')} "4";"6"};
{\ar@{=>}^<<<{\alpha} (0,3)*{};(0,-3)*{}} ;
\endxy
\]
consists of a 2-cell $\alpha : U \Rightarrow U' \in \cl{B}$ making the following diagrams commute:

\begin{equation}\label{mmtrans1}
\xy
(0,0)*+{S' U'}="1";
(0,20)*+{S' U}="2";
(20,20)*+{U S}="3";
(20,0)*+{U' S}="4";
{\ar^{\sigma} "2";"3"};
{\ar_{\sigma'} "1";"4"};
{\ar_{S' \alpha} "2";"1"};
{\ar^{\alpha S} "3";"4"};
\endxy\end{equation}

\begin{equation}\label{mmtrans2}
\xy
(0,0)*+{T' U'}="1";
(0,20)*+{T' U}="2";
(20,20)*+{U T}="3";
(20,0)*+{U' T}="4";
{\ar^{\tau} "2";"3"};
{\ar_{\tau'} "1";"4"};
{\ar_{T' \alpha} "2";"1"};
{\ar^{\alpha T} "3";"4"};
\endxy\end{equation}

\end{prop}

\noindent Note that diagram (\ref{mmtrans1}) makes $\alpha$ into a 2-cell $(U, \sigma) \Rightarrow (U', \sigma') \in \mnd{\cl{B}}$ and diagram (\ref{mmtrans2}) makes $\alpha$ into a 2-cell $(U, \tau) \Rightarrow (U', \tau') \in \mnd{\cl{B}}$.  

\begin{proof}
A 2-cell as required consists of a 2-cell in \mnd{\cl{B}}
\[
\xy
(-12,0)*+{(X, S)}="4";
(12,0)*+{(X', S')}="6";
{\ar@/^1.65pc/^{(U, \sigma)} "4";"6"};
{\ar@/_1.65pc/_{(U', \sigma')} "4";"6"};
{\ar@{=>}^<<<{\alpha} (0,3)*{};(0,-3)*{}} ;
\endxy
\]
making the following diagram commute

\begin{equation}\label{diag2}
\xy
(0,0)*+{(T', \lambda') \circ (U', \sigma'))}="1";
(0,25)*+{(T', \lambda') \circ (U, \sigma)}="2";
(40,25)*+{(U', \sigma') \circ (T, \lambda)}="3";
(40,0)*+{(U', \sigma') \circ (T', \lambda')}="4";
{\ar^{\tau} "2";"3"};
{\ar_{\tau'} "1";"4"};
{\ar_{(T',\lambda').\alpha} "2";"1"};
{\ar^{\alpha.(T,\lambda)} "3";"4"};
\endxy\end{equation}

\noindent Now the 2-cell $\alpha \in \mnd{\cl{B}}$ gives the required 2-cell $\alpha \in \cl{B}$ and diagram (\ref{mmtrans1}); diagram (\ref{diag2}) in \mnd{\cl{B}} becomes diagram (\ref{mmtrans2}) in \cl{B}.  \end{proof}


The next theorem tells us that our notation ``\mndnb{2}'' is more than just a piece of notation.  Recall that Theorem~\ref{barrwells} tells us that a distributive law of $S$ over $T$ makes $TS$ into a monad.

\begin{thm}[Street \cite{str1}]
The assignation $\mbox{{\upshape Mnd}}: \cl{B} \mapsto \mnd{\cl{B}}$ extends to a monad on the category \cat{2-Cat} of 2-categories and 2-functors, with multiplication given by
\[\begin{array}{ccc}
\mndnb{2} & \lra & \mnd{\cl{B}} \\
ST \stackrel{\lambda}{\Rightarrow} TS & \mapsto & TS.
\end{array}\]

\end{thm}

Our aim is to prove that a 0-cell of \mndnb{n} consists of a distributive series of monads $T_1, \ldots T_n$ as in Theorem~\ref{mainthm} (although we will index them as in Remark~\ref{mainthmrem}).  In order to use induction to prove this, we must characterise the whole 2-category structure.  This is the content of the next theorem.

First note that we will use the following notation
\[\big(X, \{S_i\}_{1 \leq i \leq n}, \{\lambda_{ij}\}_{1 \leq i < j \leq n}\big)\]
as a shorthand for 
\[(X, S_1, \ldots, S_n, \lambda_{12}, \lambda_{13} \ldots, \lambda_{1n}, \lambda_{23}, \lambda_{24}, \ldots, \lambda_{2n}, \ldots, \lambda_{n-1 n});\]
when the range of the indices $i,j$ is unambiguous we will simply write $\big( X, \{S_i\}, \{\lambda_{ij}\}\big)$.

\begin{thm}
Fix $n \geq 3$.  The 2-category \mndnb{n} is given as follows.  A 0-cell consists of a tuple 
\[\big(X, \{S_i\}_{1 \leq i \leq n}, \{\lambda_{ij}\}_{1 \leq i < j \leq n}\big)\]
where

\begin{itemize}

\item $X$ is a 0-cell of \cl{B},
\item each $S_i$ is a monad on $X$,  
\item each $\lambda_{ij}$ is a distributive law $S_i S_j \lra S_j S_i$, and
\item for all $i < j < k$ the monads $S_i, S_j, S_k$ satisfy the Yang-Baxter equation.
\end{itemize}

\bigskip

A 1-cell $\big( X, \{S_i\}, \{\lambda_{ij}\}\big) \lra \big( X', \{S'_i\}, \{\lambda'_{ij}\}\big)$ is a tuple $\big(T, \{\tau_i\}_{1\leq i \leq n}\big)$
where

\begin{itemize}
\item $T$ is a 1-cell $X \lra X' \in \cl{B}$,
\item each $\tau_i$ is a 2-cell 
\[S'_i T \stackrel{\tau_i}{\Rightarrow} TS_i \in \cl{B},\] 
and 
\item for all $ 1\leq i \leq n$ the following diagrams commute

\begin{equation}\label{diag3}
\xymatrix{
& T \ar[dl]_{\eta^{S'_i} T} \ar[dr]^{T\eta^{S_i}} && {S'_i}^2 T \ar[r]^<<<<{S'_i\lambda} \ar[d]_{\mu^{S'_i} T}
& S'_iTS_i \ar[r]^{\lambda S_i }
& T{S_i}^2 \ar[d]^{T \mu^{S_i}} \\
S'_iT \ar[rr]_{\lambda} && TS_i & S'_iT \ar[rr]_{\lambda}
&& TS_i}\end{equation}
i.e. each $(T,\tau_i)$ is a morphism $(X,S_i) \lra (X, S'_i) \in \mndb$, and

\item for all $i<j$ the following diagram commutes

\begin{equation}\label{diag4}
\xy
(0,12)*+{S_i' S_j' T}="1";
(12,24)*+{S_j' S_i' T}="2";
(36,24)*+{S_j' T S_i}="3";
(12,0)*+{S_i' T S_j}="4";
(36,0)*+{TS_i S_j}="5";
(48,12)*+{TS_j S_i}="6";
{\ar^{\lambda_{ij}' T} "1";"2"};
{\ar^{S_j' \tau_i} "2";"3"};
{\ar^{\tau_j S_i} "3";"6"};
{\ar_{S_i' \tau_j} "1";"4"};
{\ar_{\tau_i' S_j} "4";"5"};
{\ar_{T \lambda_{ij}} "5";"6"};
\endxy\end{equation}

\end{itemize}

\noindent The 1-cell composite
\[\xymatrix{
(X, S_i, \lambda_{ij}) \ar[r]^{(T,\tau_i)} & (X', S'_i, \lambda'_{ij}) \ar[r]^{(T'_i, \tau'_i)} & (X'', S''_i, \lambda''_{ij}) }\]

\noindent is given by $(T'T, \ T'\tau_i \circ \tau'_i T)$.  

\bigskip

\noindent A 2-cell in \mndnb{n}
\[
\xy
(-16,0)*+{(X_i, S_i, \lambda_{ij})}="4";
(16,0)*+{(X'_i, S'_i, \lambda'_{ij})}="6";
{\ar@/^1.9pc/^{(T, \tau_i)} "4";"6"};
{\ar@/_1.9pc/_{(T', \tau'_i)} "4";"6"};
{\ar@{=>}^<<<{\alpha} (0,3)*{};(0,-3)*{}} ;
\endxy
\]

\noindent consists of 

\begin{itemize}
\item a 2-cell $\alpha:T \Rightarrow T' \in \cl{B}$,  such that
\item for all $1 \leq i \leq n$ the following diagram commutes

\begin{equation}\label{diag5}
\xy
(0,0)*+{S'_i T'}="1";
(0,20)*+{S'_i T}="2";
(20,20)*+{TS_i}="3";
(20,0)*+{T'S_i}="4";
{\ar^{\tau_i} "2";"3"};
{\ar_{\tau'_i} "1";"4"};
{\ar_{S'_i \alpha} "2";"1"};
{\ar^{\alpha S_i} "3";"4"};
\endxy
\end{equation}
making $\alpha$ into a 2-cell $(T,\tau_i) \Rightarrow (T'_i, \tau'_i) \in \mnd{\cl{B}}$ for all $i$.
\end{itemize}

\noindent 2-cell composition is inherited from \cl{B}.  

\end{thm}

\noindent Note that for notational convenience in proving this theorem by induction, we have used the ``reverse'' order of indexing as in Remark~\ref{mainthmrem}.  

\begin{proof}

We write $\cl{E}_n$ for the 2-category above and prove $\cl{E}_n = \mndnb{n}$ by induction. First we prove the case $n=3$, that is, we show that

\[\mnd{\mndnb{2}} \simeq \cl{E}_3.\]

\subsubsection*{0-cells}

A 0-cell in \mnd{\mndnb{2}} consists of the following cells in \mndnb{2}:

\numarabic
\begin{enumerate}

\item a 0-cell $((X,S), (T, \lambda))$,
\item a 1-cell $((U, \sigma), \tau) : ((X,S), (T,\lambda)) \lra ((X, S), (T, \lambda))$,

\item a 2-cell $\eta: 1 \Rightarrow ((U, \sigma), \tau)$, and
\item a 2-cell $\mu: ((U, \sigma), \tau)^2 \Rightarrow ((U,\sigma), \tau)$
\end{enumerate}
satisfying the monad axioms (\ref{monadax}).  

Now (1) gives monads $S$ and $T$ and a distributive law \[\lambda:ST \Rightarrow TS.\]  The 1-cell (2) gives a 1-cell $U:X \lra X$ that is made into a monad by (3) and (4).  The 1-cell (2) also gives a 2-cell \[\sigma: SU \Rightarrow US\] satisfying the first pair of distributive law axioms, diagrams (\ref{dl1}) governing interaction with the monad structure of $S$; the other pair of axioms, diagrams (\ref{dl2}), come from (3) and (4).  Further, (2) gives a 2-cell \[\tau: TU \Rightarrow UT\] satisfying the Yang-Baxter equation for $S,T,U$.  Axioms (\ref{mmmorph2}) for a 1-cell ensure that $\tau$ interacts properly with the monad structure for $T$; (3) and (4) above ensure that $\tau$ interacts properly with the monad structure for $U$, hence is a distributive law.  This gives the result, where we have written $S,T,U$ for $S_1, S_2, S_3$, and similarly for $\lambda, \sigma, \tau$.

\subsubsection*{1-cells}

A 1-cell in \mnd{\mndnb{2}} 
\[ (X, \{S, T, U\}, \{\lambda, \sigma, \tau\}) \lra (X', \{S', T', U'\}, \{\lambda', \sigma', \tau'\}) \]
consists of the following cells in \mndnb{2}

\numarabic
\begin{enumerate}

\item a 1-cell $((V, \theta), \phi) : ((X,S), (T,\lambda)) \lra ((X',S'), (T',\lambda'))$, and
\item a 2-cell $\rho: ((U',\sigma'), \tau') \circ ((V, \theta), \phi) \lra ((V,\theta), \phi) \circ ((U, \sigma), \tau)$
\end{enumerate}

\noindent such that the following diagrams commute (where we now omit the sub-parentheses for convenience)
\[\xy
(25,15)*+{(V,\theta,\phi)}="1";
(0,0)*+{(U', \sigma', \tau') \circ (V, \theta, \phi)}="2";
(50,0)*+{(V,\theta,\phi) \circ (U, \sigma, \tau)}="3";
{\ar_{\eta.1} "1";"2"};
{\ar^{1.\eta} "1";"3"};
{\ar_{\rho} "2";"3"};
\endxy\]


\begin{equation}\label{diag6}
\xy
(0,20)*+{(U',\sigma',\tau')^2 \circ (V,\theta,\phi)}="1";
(52,20)*+{(U',\sigma',\tau')\circ (V,\theta,\phi)\circ (U, \sigma, \tau)}="2";
(104,20)*+{(V, \theta, \phi)\circ (U, \sigma, \tau)^2}="3";
(0,0)*+{(U', \sigma', \tau') \circ (V, \theta, \phi)}="4";
(104,0)*+{(V,\theta,\phi)\circ (U, \sigma, \tau)}="5";
{\ar^<<<<<<{1.\rho} "1";"2"};
{\ar^>>>>>>>{\rho.1} "2";"3"};
{\ar^{\mu.1} "1";"4"};
{\ar^{1.\mu} "3";"5"};
{\ar_{\rho} "4";"5"};
\endxy\end{equation}

\vspace{2em}\noindent Now (1) gives 1-cells in \mnd{\cl{B}} as follows
\[(V,\theta): (X,S) \lra (X',S')\] 
and 
\[(V,\phi):(X,T) \lra (X',T')\] 
such that the hexagon~(\ref{mmmorph3}) commutes for $S,T,V$.  The 2-cell (2) gives a 2-cell 
\[\rho: U' V \lra VU \in \cl{B}\] 
which is made into a monad map 
\[(V,\rho):(X,U) \lra (X',U')\] 
by diagrams (\ref{diag6}).  Further, (2) has two diagrams (\ref{mmtrans1}) and (\ref{mmtrans2}); diagram (\ref{mmtrans1}) becomes the hexagon (\ref{mmmorph3}) for $S,U,V$ and diagram (\ref{mmtrans2}) becomes the hexagon for $T,U,V$.  It is straightforward to check the formula for composition, so this completes the result for 1-cells.

\subsubsection*{2-cells}

A 2-cell in \mnd{\mndnb{2}} consists of

\begin{enumerate}
\item a 2-cell $\alpha: (V,\theta,\phi) \lra (V', \theta',\phi') \in \mndnb{2}$, such that
\item the following diagram commutes

%
\[\xy
(0,25)*+{(U',\sigma',\tau') \circ (V,\theta,\phi)}="1";
(45,25)*+{(V,\theta,\phi) \circ (U,\sigma,\tau)}="2";
(0,0)*+{(U', \sigma', \tau') \circ (V', \theta', \phi') }="3";
(45,0)*+{(V', \theta', \phi') \circ (U, \sigma, \tau)}="4";
{\ar^{\rho} "1";"2"};
{\ar^{\alpha.(U, \sigma, \tau)} "2";"4"};
{\ar_{(U', \sigma', \tau').\alpha} "1";"3"};
{\ar_{\rho'} "3";"4"};
\endxy\]
\end{enumerate}

\noindent Now (1) gives $\alpha$ as a 2-cell 
\[(V,\theta) \Rightarrow (V',\theta') \in \mnd{\cl{B}} \]
and also as a 2-cell 
\[(V,\phi) \Rightarrow (V',\phi') \in \mnd{\cl{B}}; \] 
diagram (2) makes $\alpha$ into a 2-cell 
\[(V,\rho) \Rightarrow (V',\rho') \in \mndb.\]  
This completes the result for 2-cells and thus the case for $n=3$.  

We now prove the case for $n$, that is, that $\mnd{\cl{E}_{n-1}}=\cl{E}_n$.  

\subsubsection*{0-cells of \mnd{\cl{E}_{n-1}}}

A 0-cell of \mnd{\cl{E}_{n-1}} consists of:

\begin{enumerate}

\item a 0-cell $(X, \{S_i\}_{1\leq i \leq n-1}, \{\lambda_{ij}\}_{1\leq i < j \leq n-1})$ of $\cl{E}_{n-1}$,
\item a 1-cell $(T,\{\tau_i\}_{1\leq i \leq n-1}): (X, \{S_i\}, \{\lambda_{ij}\}) \lra (X, \{S_i\}, \{\lambda_{ij}\})$, and
\item 2-cells $\eta: (1,1) \Rightarrow (T,\{\tau_i\})$ and $\mu: (T,\{\tau_i\})^2 \Rightarrow (T,\{\tau_i\})$
\end{enumerate}
satisfying the usual monad axioms (\ref{monadax}).

Note that {\em a priori} our indices only run from 1 to $n-1$, so it remains to define a monad $S_n$, distributive laws \[\lambda_{in}: S_i S_n \Rightarrow S_n S_i\] for all $1 \leq i \leq n-1$, and check the Yang-Baxter equations for all triples of monads $S_i, S_j, S_n$.  

Now (2) certainly gives a 1-cell $T:X \lra X \in \cl{B}$ which is made into a monad by the underlying 2-cells of (3) and the monad axioms.  So we put $S_n = T$. The 1-cell (2) also gives for each $i \leq n-1$ a morphism
\[(T,\tau_i): (X, S_i) \lra (X,S_i) \in \mndb \]

\noindent with 2-cell component \[\tau_i: S_i T \lra TS_i,\] so we set $\lambda_{in} = \tau_i$ for each $i<n$.  The axioms for a monad map (\ref{mmorph}) give the interaction of each $\lambda_{in}$ with the monad structure of $S_i$, and the 2-cell axioms~(\ref{mtrans}) for $\eta$ and $\mu$ give the interaction of each $\lambda_{in}$ with the monad structure for $S_n$.  So we have all the required distributive laws $\lambda_{ij}$.  Furthermore (2) gives for all $1\leq i<j \leq n-1$ the hexagon~(\ref{diag4}) for $S_i, S_j, T$, i.e. the Yang-Baxter equation.  So we have all the required Yang-Baxter equations.

\subsubsection*{1-cells of \mnd{\cl{E}_{n-1}}}

A 1-cell of \mnd{\cl{E}_{n-1}} 
\[ \big( (X, \{S_i\}, \{\lambda_{ij}\}),\ (S_n, \{\lambda_{in}\}) \big) \lra \big( (X', \{S'_i\}, \{\lambda'_{ij}\}),\ (S'_n, \{\lambda'_{in}\}) \big)\]

\noindent consists of the following cells in $\cl{E}_{n-1}$:

\begin{enumerate}

\item a 1-cell $\big(U, \{\sigma_i\}_{1\leq i \leq n-1}\big) : \big(X, \{S_i\}, \{\lambda_{ij}\}\big) \lra  \big(X', \{S'_i\}, \{\lambda'_{ij}\}\big)$, and

\item a 2-cell $\alpha: \big(S'_n, \{ \lambda'_{in}\}\big) \circ \big(U, \{\sigma_i\}\big) \Rightarrow 
\big(U, \{\sigma_i\}\big) \circ \big(S_n, \{ \lambda_{in}\}\big) $
\end{enumerate}

\noindent satisfying diagrams~(\ref{mmorph}).  As in the case of the 0-cells, it remains to define a 2-cell \[\sigma_n : S'_nU \Rightarrow US_n \in \cl{B},\] and check diagrams~(\ref{diag3}) and~(\ref{diag4}).  Now (2) gives a 2-cell
\[\alpha : \big(S'_nU, \{U\lambda'_{in} \circ S'_n \sigma_i\}\big) \Rightarrow \big(US_n, \{S_n \sigma_i \circ U\lambda_{in}\}\big) \in \cl{E}_{n-1}\]

\noindent thus a 2-cell \[\alpha : S'_nU \Rightarrow US_n \in \cl{B},\] and diagrams~(\ref{mmorph}) make $(U,\alpha)$ into a monad map \[(U, \alpha):(X,S_n) \lra (X', S'_n).\]  So we set $\sigma_n = \alpha$.  The axioms~(\ref{diag5}) making $\alpha$ a 2-cell of $\cl{E}_{n-1}$ give the hexagon~(\ref{diag4}) for $S_i, S_n, U$ for all $1\leq i <n$; the other hexagons come from (1).  Finally 1-cell composition is given by

\begin{eqnarray*}
\big((U', \{\sigma'_i\}), \alpha'\big) \circ \big((U, \{\sigma_i\}), \alpha \big) & = & \big((U', \{\sigma'_i\}) \circ (U, \{\sigma_i\}) , \ (U', \{\sigma'_i\}).\alpha \circ \alpha. (U, \{\sigma_i\}) \big) \\
& = & \big( (U'U, \ U'\{\sigma_i\} \circ U \{\sigma'_i\}), \ U'\alpha \circ \alpha U\big)
\end{eqnarray*}
as required.

\subsubsection*{2-cells of \mnd{\cl{E}_{n-1}}}

A 2-cell 
\[
\xy
(-20,0)*+{\big(X, \{S_i\}, \{\lambda_{ij}\}\big)}="1";
(20,0)*+{\big(X', \{S'_i\}, \{\lambda'_{ij}\}\big)}="2";
{\ar@/^1.9pc/^{(U, \{\sigma_i\})} "1";"2"};
{\ar@/_1.9pc/_{(U', \{\sigma'_i\})} "1";"2"};
{\ar@{=>}^<<<{\rho} (0,3)*{};(0,-3)*{}} ;
\endxy
\]
in \mnd{\cl{E}_{n-1}} consists of a 2-cell 
\[\rho : (U, \{\sigma_i\}_{1\leq i \leq n-1}) \Rightarrow (U', \{\sigma'_i\}_{1\leq i \leq n-1}) \in \cl{E}_{n-1}\]
such that the following diagram commutes:
\[\xy
(0,0)*+{(S'_n, \{\lambda_{in}\}) \circ (U', \{\sigma'_i\}) }="1";
(0,25)*+{(S'_n, \{\lambda'_{in}\}) \circ (U, \{\sigma_i\})}="2";
(50,25)*+{(U,\{\sigma_i\}) \circ (S_n, \{\lambda_{in}\})}="3";
(50,0)*+{(U', \{\sigma'_i\}) \circ (S_n, \{\lambda_{in}\}) }="4";
{\ar^{\sigma_n} "2";"3"};
{\ar^{\rho.(S_n, \{\lambda_{in}\})} "3";"4"};
{\ar_{(S'_n, \{\lambda'_{in}\}).\rho} "2";"1"};
{\ar_{\sigma'_n} "1";"4"};
\endxy\]

\noindent Now $\rho$ being a 2-cell of $\cl{E}_{n-1}$ tells us that for all $i<n$, $\rho$ is a 2-cell \[\rho:(U,\sigma_i) \Rightarrow (U', \sigma'_i) \in \mndb;\] for the case $i=n$ the commutative diagram gives us that $\rho$ is a 2-cell \[\rho:(U,\sigma_n) \Rightarrow (U', \sigma'_n) \in \mndb\] giving us the desired result.  

\end{proof}

\subsection{Simplicial resolution of a monad}\label{simplicialres}

In this section we briefly discuss the simplicial resolution of a monad and how applying this to the monad Mnd sheds light on the results of the previous section.

Recall that given any monad $T$ we can construct its simplicial resolution:
\[\xy
(120,0)*+{T}="1"; 
(100,0)*+{T^2}="2"; 
(80,0)*+{T^3}="3"; 
(60,0)*+{T^4}="4"; 
(40,0)*+{\cdots}="5"; 
(20,0)*+{}="6"; 
(0,0)*+{T^n}="7";
(-5,0)*+{\cdots}="8";
{\ar^{\mu} "2"; "1"};
{\ar@<1ex>^{T\mu} "3"; "2"};
{\ar@<-1ex>_{\mu T} "3"; "2"};
{\ar@<3ex>^{T^2\mu} "4"; "3"};
{\ar@<0ex>^{T\mu T} "4"; "3"};
{\ar@<-3ex>_{\mu T^2} "4"; "3"};
{\ar@<4.5ex>^{} "5"; "4"};
{\ar@<1.5ex>^{} "5"; "4"};
{\ar@<-1.5ex>^{} "5"; "4"};
{\ar@<-4.5ex>^{} "5"; "4"};
{\ar@<-7ex>^{} "7"; "6"};
{\ar@<-5ex>^{} "7"; "6"};
{\ar@<-3ex>^{} "7"; "6"};
{\ar@<-1ex>^{} "7"; "6"};
{\ar@<1ex>^{} "7"; "6"};
{\ar@<3ex>^{\vdots} "7"; "6"};
{\ar@<9ex>^{} "7"; "6"};
\endxy\]

\noindent with various commuting conditions ensuring, among other things, that the diagram yields a unique morphism from $T^n$ to $T$.  Note that we have only drawn the multiplications (face maps) in this diagram; there are also degeneracies corresponding to applications of the unit for the monad. Applying this construction to the monad Mnd, we see that the unique morphism from $\mbox{Mnd}^n$ to Mnd gives us the unique composite monad $T_1 T_2 \cdots T_n$; the maps to $\mbox{Mnd}^2$ give us the distributive laws \ $T_{i+1} T_{i+2} \cdots T_n$\  over \ $T_1 T_2 \cdots T_i$.  



Furthermore, we can use the simplicial resolution of monads to express monads in \cl{B} as 2-functor from a certain 2-category $\Delta$ to \cl{B}.  

Let $\Delta$ be the ``free-living monad'' 2-category of ordinals.  $\Delta$ is more commonly thought of as a \emph{category} whose objects are the natural numbers (including 0), but it has a monoidal structure given by addition; thus it can be considered as a bicategory with only one 0-cell. Then a monad in a 2-category \cl{B} can be expressed as a (strict) functor $\Delta \lra \cl{B}$. The image of the single 0-cell of $\Delta$ picks out an underlying 0-cell $X$ of \cl{B} and the rest of $\Delta$ picks out a monad on $X$ by specifying its entire simplicial resolution.

Furthermore,  a monad map is a lax transformation between functors, and a monad transformation is a modification.  Recall  \cite{gra1, gra2} that this combination of strictness and laxness gives us a closed structure with respect to the lax Gray tensor product as follows.

Given 2-categories $\cl{A}$ and $\cl{B}$, write $[\cl{A},\cl{B}]$ for the 2-category whose 0-cells are strict functors $\cl{A} \lra \cl{B}$, 1-cells are lax transformations and 2-cells are modifications.  Write $\cl{A} \otimes \cl{B}$ for the lax Gray tensor product of $\cl{A}$ and $\cl{B}$, and $\cat{Gray}_{\mbox{lax}}$ for the monoidal category of 2-categories and 2-functors with monoidal structure given by the the lax Gray tensor product.  Then $\cat{Gray}_{\mbox{lax}}$ is closed with internal hom given by $[\cl{A},\cl{B}]$.

Thus we have:
\[\begin{array}{lclcl}
\mndb & = & [\Delta, \cl{B}] \\[2ex]
\mndnb{2} & = & \big[\Delta,\ [\Delta, \cl{B}]\big] &=& [\Delta \otimes \Delta, \cl{B}] \\
& \vdots \\
\mndnb{n} &=& [\Delta \otimes \cdots \otimes \Delta, \cl{B}] &=& [\Delta^{\otimes n}, \cl{B}]
\end{array}\]

\noindent We now sketch the correspondence
\[[\Delta \otimes \Delta, \cl{B}]  \lra  \mndnb{2}.\]
Let us write the 1-cells of $\Delta$ as $1, e, e^2, e^3, \dots $.  Then $\Delta \otimes \Delta$ has 1-cells generated by $(1,e)$ and $(e,1)$.  Now whereas in $\Delta \times \Delta$ we have the relation
\[(1,e) \circ (e,1) = (e,1) \circ (1,e),\]
in $\Delta \otimes \Delta$ we instead have a 2-cell generator
\[\xy
(0,0)*+{.}="1";
(0,20)*+{.}="2";
(20,20)*+{.}="3";
(20,0)*+{.}="4";
{\ar^{(1,e)} "2";"3"};
{\ar_{(1,e)} "1";"4"};
{\ar_{(e,1)} "2";"1"};
{\ar^{(e,1)} "3";"4"};
{\ar@{=>}^{\phi} (12,12);(8,8)};
\endxy\]

\noindent We show how a functor \[\alpha: \Delta \otimes \Delta \lra \cl{B}\] corresponds to an object of \mndnb{2}, that is, monads $S$ and $T$ and a distributive law $\lambda: ST \Rightarrow TS$.  

First recall that a functor $\theta: \Delta \lra \cl{B}$ gives us a monad in \cl{B} by giving us the entire simplicial resolution of a monad, thus the image of $e$ gives the functor part of the monad.  So for the case above we can set $S = \alpha (1,e)$ and $T = \alpha(e,1)$ and these are automatically monads.  For the distributive law recall that in $\Delta \otimes \Delta$ we have the 2-cell $\phi$ above.  Now $\alpha$ applied to the upper right leg of the square gives $ST$, and applied to the lower left leg it gives $TS$.  So we have \[\alpha(\phi) : ST \Rightarrow TS\] and this can be shown to be a distributive law.  Furthermore, to find the composite monad $TS$ we use a ``diagonal'' functor:
\[\begin{array}{ccc}
\Delta & \map{d} & \Delta \otimes \Delta\\
f & \mapsto & (1,f) \circ (f,1) 
\end{array}\]

\noindent Then given any $\alpha : \Delta \otimes \Delta \Rightarrow \cl{B}$ corresponding to $(S,T,\lambda)$ we get the lax functor
\[\Delta \map{d} \Delta \otimes \Delta \map{\alpha} \cl{B}\]
corresponding to the monad $TS$.  

Finally note that in the definition of $\Delta \otimes \Delta \otimes \Delta$ by generators and relations, a Yang-Baxter equation is seen to arise from the relations; this corresponds to the Yang-Baxter equation we have seen in \mndnb{n}.

\section{Interchange for $n$-categories}\label{interchange}

In this section we discuss our motivating example, the free strict $n$-category monad on $n$-dimensional globular sets.  In an $n$-category, interchange laws govern the interaction between different types of composition.  These different types of composition can be expressed using monads, and the main result of this section is that the interchange laws define distributive laws between those monads, giving a distributive series of monads.  Using the theory of iterated distributive laws, the resulting composite monad is the standard ``free strict $n$-category'' monad induced by the adjunction:
\[\xy
(0,0)*+{\cat{$n$-GSet}}="1";
(30,0)*+{\cat{Str-$n$-Cat}}="2";
(14,0)*{{\scriptstyle \perp}};
{\ar@<1ex>_{} "1";"2"};
{\ar@<1ex>^{} "2";"1"};
\endxy\]

Throughout this section we will omit the word ``strict'' and understand all our $n$-categories to be strict.  In fact the key for us is that \emph{interchange} is strict; this theory could in principle be used for notions of $n$-category that are weaker, as long as interchange is still strict.  An example of this is Trimble's definition \cite{tri1, che16}.

\subsection{Composition in $n$-categories}

The underlying data for an $n$-category is an $n$-globular set, that is, a diagram of sets and functions
\[\xy
(120,0)*+{A(0)}="0"; 
(100,0)*+{A(1)}="1"; 
(80,0)*+{A(2)}="2"; 
(63,0)*+{}="3"; 
(50,0)*+{}="4"; 
(30,0)*+{A(n-1)}="5"; 
(5,0)*+{A(n)}="6";
(57,0)*{\cdots};
{\ar@<0.7ex>^{s} "1"; "0"};
{\ar@<-0.7ex>_{t} "1"; "0"};
{\ar@<0.7ex>^{s} "2"; "1"};
{\ar@<-0.7ex>_{t} "2"; "1"};
{\ar@<0.7ex>^<<<<<<{s} "3"; "2"};
{\ar@<-0.7ex>_<<<<<<{t} "3"; "2"};
{\ar@<0.7ex>^>>>>>>>{s} "5"; "4"};
{\ar@<-0.7ex>_>>>>>>>{t} "5"; "4"};
{\ar@<0.7ex>^<<<<<<<{s} "6"; "5"};
{\ar@<-0.7ex>_<<<<<<<{t} "6"; "5"};
\endxy\]
such that $ss=st, ts=tt$. The elements of each $A(m)$ are generally referred to as ``$m$-cells'', and the functions $s$ and $t$ give the ``source'' and ``target'' $m$-cells of an $(m+1)$-cell, also generally known as the boundary or bounding cells.  Then $n$-globular sets form a category \cat{$n$-GSet} with the obvious morphisms; note that \cat{$n$-GSet} can be expressed as the category of presheaves in the obvious way.   

An $n$-category should be an $n$-globular set with, for all $0 \leq m \leq n-1$, composition along bounding $m$-cells, which we will call $m$-composition and denote by $\circ_m$.  For example 2-categories have:

\begin{itemize}

\item 0-composition = horizontal composition, usually denoted $*$, or in diagrams
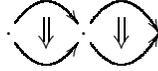
\[\xy
(0,0)*+{.}="1";
(10,0)*+{.}="2";
(20,0)*{}="3";
{\ar@/^1pc/^{} "1";"2"};
{\ar@/_1pc/_{} "1";"2"};
{\ar@{=>}^{} (5,2)*{};(5,-2)*{}} ;
{\ar@/^1pc/^{} "2";"3"};
{\ar@/_1pc/_{} "2";"3"};
{\ar@{=>}^{} (15,2)*{};(15,-2)*{}} ;
\endxy
\]

\item 1-composition = vertical composition, usually denoted $\circ$, or in diagrams
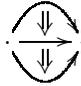
\[\xy
(0,0)*+{.}="1";
(10,0)*+{.}="2";
{\ar@/^1.25pc/^{} "1";"2"};
{\ar@/_1.25pc/_{} "1";"2"};
{\ar "1";"2"};
{\ar@{=>}^{} (5,4)*{};(5,1)*{}} ;
{\ar@{=>}^{} (5,-1)*{};(5,-4)*{}} ;
\endxy\]

\end{itemize}
and each is strictly unital and associative.  In a 2-category we also have the interchange law
\[(a * b) \circ (c*d) = (a \circ c) * (b \circ d)\]
or in diagrams
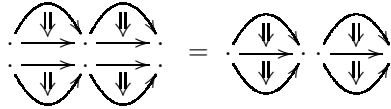
\[\xy
(30,0)*{ 	
\xy
(0,0)*{
\xy
(0,0)*+{.}="1";
(10,0)*+{.}="2";
{\ar@/^1.25pc/^{} "1";"2"};
{\ar@/_1.25pc/_{} "1";"2"};
{\ar "1";"2"};
{\ar@{=>}^{} (5,4)*{};(5,1)*{}} ;
{\ar@{=>}^{} (5,-1)*{};(5,-4)*{}} ;
\endxy};
(12,0)*{
\xy
(0,0)*+{.}="1";
(10,0)*+{.}="2";
{\ar@/^1.25pc/^{} "1";"2"};
{\ar@/_1.25pc/_{} "1";"2"};
{\ar "1";"2"};
{\ar@{=>}^{} (5,4)*{};(5,1)*{}} ;
{\ar@{=>}^{} (5,-1)*{};(5,-4)*{}} ;
\endxy}
\endxy};
(15,0)*{=};
(0,0)*{ 
\xy
(0,0)*{
\xy
(0,0)*+{.}="1";
(10,0)*+{.}="2";
(20,0)*+{.}="3";
{\ar@/^1.25pc/^{} "1";"2"};
{\ar "1";"2"};
{\ar@{=>}^{} (5,4)*{};(5,1)*{}} ;
{\ar@/^1.25pc/^{} "2";"3"};
{\ar "2";"3"};
{\ar@{=>}^{} (15,4)*{};(15,1)*{}} ;
\endxy};
(0,-8)*{
\xy
(0,0)*+{.}="1";
(10,0)*+{.}="2";
(20,0)*+{.}="3";
{\ar@/_1.25pc/^{} "1";"2"};
{\ar "1";"2"};
{\ar@{=>}^{} (5,-1)*{};(5,-4)*{}} ;
{\ar@/_1.25pc/^{} "2";"3"};
{\ar "2";"3"};
{\ar@{=>}^{} (15,-1)*{};(15,-4)*{}} ;
\endxy};
\endxy};
\endxy\]
In effect this, together with associativity and unit laws, ensures that any given diagram of composable cells has a unique composite.  For $m$-cells in an $n$-category there are $m$ different kinds of composition, along bounding $i$-cells for all $0 \leq i \leq m-1$, and an interchange law for all pairs $i,j$ with $0 \leq i < j \leq m-1$
\[(a \circ_j b) \circ_i (c \circ_j d) = (a \circ_i c) \circ_j (b \circ_i d)\]
ensuring that any diagram of composable cells (perhaps including more than two types of composition) has a unique composite.

\subsection{Monads for $i$-composition}

We construct, for each $0 \leq i < n$ a monad $T_i$ on \cat{$n$-GSet} which constructs $i$-composites freely (leaving $k$-cells alone for $k \leq i$).  Each of these monads is a completely straightforward generalisation of the ordinary free category monad on graphs.  We give the details here simply in order to be able to show that the composite monad $T_0 T_2 \cdots T_{n-1}$ resulting from the distributive series of monads in question, is really the free strict $n$-category monad.  We draw on abstract results from Appendix F of \cite{lei8}, but writing down the definitions directly is not hard.

The construction of the monad for $i$-composition proceeds in the following steps:

\numarabic
\begin{enumerate}
\item Construct free category monad on \cat{1-GSet}.
\item Use the enriched version to construct a monad for ``free 0-composition''  on \cat{$(n-i)$-GSet}.
\item Shift the dimensions up $i$ times by inserting lower dimensions trivially, which turns this into the monad for ``free $i$-composition'' on \cat{$n$-GSet}.
\end{enumerate}

First we recall the free enriched category monad as described in \cite{lei8}, which acts on the category of \cl{V}-graphs.

\begin{Def}

Given a category \cl{V}, a {\em \cl{V}-graph} $A$ is given by 

\begin{itemize}
\item a set $A_0$ of objects, and 
\item for every pair of objects $a,a'$, an object $A(a,a') \in \cl{V}$.
\end{itemize}
A morphism $F: A \lra B$ of \cl{V}-graphs is given by

\begin{itemize}
\item a function $F:A_0 \lra B_0$, and 
\item for every pair of objects $a,a'$, a morphism $A(a,a') \lra B(Fa, Fa') \in \cl{V}$.
\end{itemize}
\cl{V}-graphs and their morphisms form a category \cat{\cl{V}-Gph}.  
\end{Def}
Note that
\[\begin{array}{ccc}
\cat{Set-Gph} & = & \cat{1-GSet}\\
\cat{Gph-Gph} & = & \cat{2-GSet}\\
\vdots\\
\cat{($n$-GSet)-Gph} & = & \cat{$(n+1)$-GSet}
\end{array}\]
We will also write \cat{$n$-Gph} for \cat{$n$-GSet}, so $\cat{($n$-Gph)-Gph}=\cat{$(n+1)$-Gph}$.

If \cl{V} is monoidal we can construct categories enriched in \cl{V}, but to make the free \cl{V}-category construction we need \cl{V} to be suitably well-behaved.  If \cl{V} is a presheaf category it is certainly well enough behaved \cite{lei8}, thus \cat{$n$-GSet} is suitable.  Recall that a monad is called cartesian if it preserves pullbacks and the naturality squares for $\eta$ and $\mu$ are all pullbacks.  

The following theorem gives us the enriched version of the free category monad.

\begin{thm}[Leinster \cite{lei8}]\label{freevcat}

If \cl{V} is a presheaf category then the forgetful functor
\[\cat{\cl{V}-Cat} \lra \cat{\cl{V}-Gph}\]
is monadic.  The induced monad is the ``free \cl{V}-category monad'' $\cat{fc}_{\cl{V}}$ and is cartesian.

\end{thm}

The following corollary is the example we need, giving us the monad on \cat{$n$-Gph} for ``free 0-composition''.  We use $\cl{V} = \cat{$(n-1)$-Gph}$. 

\begin{cor}\label{free}
For all $n \geq 1$ we have a monadic adjunction
\[\xy
(0,0)*+{\cat{($(n-1)$-Gph)-Cat}}="1";
(60,0)*+{\cat{($(n-1)$-Gph)-Gph}=\cat{$n$-Gph}}="2";
(26,0)*{{\scriptstyle \perp}};
{\ar@<1ex>_{} "1";"2"};
{\ar@<1ex>^{} "2";"1"};
\endxy\]
The induced monad $T$ constructs 0-composites freely:
\[TA(m) = \coprod_{\begin{array}{c} {\scriptstyle k\geq 0} \\ {\scriptstyle a_0, \ldots, a_k \in A(0)} \end{array}} A(a_{k-1}, a_k) \times \cdots \times A(a_0, a_1) \]
\end{cor}
\begin{proof}
Put $\cl{V}=\cat{$(n-1)$-Gph}$ in Theorem~\ref{freevcat}.  Then $T = \cat{fc}_{\cat{{\small ${\scriptstyle (n-1)}$-Gph}}}$ and the formula is exactly the formula given in \cite{lei8}.
\end{proof}

Note that this formula produces $k$-length strings of 0-composable cells.  It is a coproduct over $k$ 
of $k$-fold wide pullbacks as below:

%
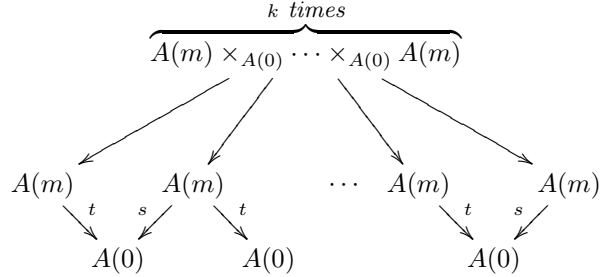
\[\xy
(0,10)*+{A(m)}="1";
(10,0)*+{A(0)}="2";
(20,10)*+{A(m)}="3";
(30,0)*+{A(0)}="4";
(50,10)*+{A(m)}="5";
(60,0)*+{A(0)}="6";
(70,10)*+{A(m)}="7";
(35,30)*+{\overbrace{A(m) \times_{A(0)} \cdots \times_{A(0)} A(m)}^{k \mbox{ {\small\it times}}}
}="8";
(40,10)*+{\cdots};
{\ar^{t} "1";"2"};
{\ar_{s} "3";"2"};
{\ar^{t} "3";"4"};
{\ar^{t} "5";"6"};
{\ar_{s} "7";"6"};
{\ar^{} "8";"1"};
{\ar_{} "8";"3"};
{\ar^{} "8";"5"};
{\ar^{} "8";"7"};
\endxy\]
where $s$ and $t$ denote the composites along the top and bottom of 
\[\xy
(120,0)*+{A(0)}="0"; 
(100,0)*+{A(1)}="1"; 
(80,0)*+{A(2)}="2"; 
(63,0)*+{}="3"; 
(50,0)*+{}="4"; 
(30,0)*+{A(m-1)}="5"; 
(5,0)*+{A(m)}="6";
(57,0)*{\cdots};
{\ar@<0.7ex>^{s} "1"; "0"};
{\ar@<-0.7ex>_{t} "1"; "0"};
{\ar@<0.7ex>^{s} "2"; "1"};
{\ar@<-0.7ex>_{t} "2"; "1"};
{\ar@<0.7ex>^<<<<<<{s} "3"; "2"};
{\ar@<-0.7ex>_<<<<<<{t} "3"; "2"};
{\ar@<0.7ex>^>>>>>>>{s} "5"; "4"};
{\ar@<-0.7ex>_>>>>>>>{t} "5"; "4"};
{\ar@<0.7ex>^<<<<<<<{s} "6"; "5"};
{\ar@<-0.7ex>_<<<<<<<{t} "6"; "5"};
\endxy\]

Note that, rather than using the abstract theory, we could simply define the monad by the formula given above and prove the later results by checking the formulae directly.  

Now in order to make free $i$-composites and not just free 0-composites we just need to ``shift'' the monad up $i$ dimensions.  The following construction shifts the monad up 1 dimension.  For any functor $F: \cl{V} \lra \cl{W}$ we get a functor $F_*: \cat{\cl{V}-Gph} \lra \cat{\cl{W}-Gph}$ as follows.  Given a $\cl{V}$-graph $A$, the graph $F_* A$ is defined by:

\begin{itemize}
\item $(F_*A)_0 = A_0$
\item $(F_*A)(a,a') = F(A(a,a'))$
\end{itemize}
and we extend this to morphisms in the obvious way.

In fact we have a 2-functor $\Cat \lra \Cat$ sending \cl{V} to $\cat{\cl{V}-Gph}$, $F$ to $F_*$ and a natural transformation $\alpha$ to a natural transformation $\alpha_*$ with components $F_* A \lra G_* A$ given by

\begin{itemize}
\item on objects the identity, which makes sense since $(F_* A)_0 = A_0 = (G_* A)_0$
\item on hom-objects \[\alpha_{A(a,a')} : F(A(a,a')) = (F_* A)(a,a') \lra (G_* A)(a,a') = G(A(a,a')).\]
\end{itemize}

\noindent This will later enable us to apply the $(-)_*$ construction to distributive laws.
The following proposition tells us that the $(-)_*$ construction preserves monadic adjunctions.

\begin{prop}[Leinster \cite{lei8}]\label{suspend}
A monadic adjunction 
\[\xy
(0,0)*+{\cl{V}}="1";
(20,0)*+{\cl{W}}="2";
(10,0)*{{\scriptstyle \perp}};
{\ar@<1ex>^{F} "1";"2"};
{\ar@<1ex>^{U} "2";"1"};
\endxy\]
induces a monadic adjunction
\[\xy
(0,0)*+{\cat{\cl{V}-Gph}}="1";
(30,0)*+{\cat{\cl{W}-Gph}}="2";
(15,0)*{{\scriptstyle \perp}};
{\ar@<1ex>^{F_*} "1";"2"};
{\ar@<1ex>^{U_*} "2";"1"};
\endxy\]
\end{prop}
\noindent Writing $T$ for the original monad $FU$, the induced monad is given by $F_* U_* = T_*$.  

\begin{eg}\label{ex1} {\bfseries 2-categories.}

\noindent We put

\begin{itemize}

\item $ \cl{V} = \cat{Gph}$,
 
\item $\cl{W} = \cat{Cat}$,

\item $F = $ the free category functor, and

\item $U = $ the usual forgetful functor.

\end{itemize}

\noindent Then Proposition~\ref{suspend} gives us an adjunction
\[\xy
(0,0)*+{\cat{2-Gph}}="1";
(30,0)*+{\cat{Cat-Gph}}="2";
(15,0)*{{\scriptstyle \perp}};
{\ar@<1ex>^{F_*} "1";"2"};
{\ar@<1ex>^{U_*} "2";"1"};
\endxy.\]
Here the functor $F_*$ sends the graph
\[\xy
(120,0)*+{A(0)}="0"; 
(100,0)*+{A(1)}="1"; 
(80,0)*+{A(2)}="2"; 
{\ar@<0.7ex>^{s} "1"; "0"};
{\ar@<-0.7ex>_{t} "1"; "0"};
{\ar@<0.7ex>^{s} "2"; "1"};
{\ar@<-0.7ex>_{t} "2"; "1"};
\endxy\]
to the cat-graph with underlying 2-graph
\[\xy
(120,0)*+{A(0)}="0"; 
(100,0)*+{A(1)}="1"; 
(80,0)*+{FA(2)}="2"; 
{\ar@<0.7ex>^{s} "1"; "0"};
{\ar@<-0.7ex>_{t} "1"; "0"};
{\ar@<0.7ex>^{s} "2"; "1"};
{\ar@<-0.7ex>_{t} "2"; "1"};
\endxy\]
where by abuse of notation we have written
\[\xy
(120,0)*+{A(1)}="0"; 
(100,0)*+{FA(2)}="1"; 
{\ar@<0.7ex>^{s} "1"; "0"};
{\ar@<-0.7ex>_{t} "1"; "0"};
\endxy\]
to denote the graph of the free category on
\[\xy
(120,0)*+{A(1)}="0"; 
(100,0)*+{A(2)}="1"; 
{\ar@<0.7ex>^{s} "1"; "0"};
{\ar@<-0.7ex>_{t} "1"; "0"};
\endxy\]
so in effect we are forming 1-composites of 2-cells freely.  This naturally has the structure of a cat-graph.  The monad $T_*$ induced by this adjunction is the free 2-category monad, and $T_*\mbox{-Alg} = \cat{2-Cat}$.

\end{eg}

We now combine Corollary~\ref{free} and Proposition~\ref{suspend} to construct the monads for $i$-composition that we require.

\begin{prop}
Let $n \geq 1$.  Then for all $0 \leq i \leq n-1$ we have a monadic adjunction
\[\xy
(0,-3)*+{\mbox{$\underbrace{\cat{Gph-Gph-$\cdots$-Gph}}_{(n-i-1) \mbox{ {\small\it times}}}$-\cat{Cat}-$\underbrace{\cat{Gph-$\cdots$-Gph}}_{i \mbox{ {\small\it times}}}$}
};
(33,0)="1";
(60,0)*+{\cat{$n$-Gph}}="2";
(43,0)*{{\scriptstyle \perp}};
{\ar@<1ex>_{} "1";"2"};
{\ar@<1ex>^{} "2";"1"};
\endxy\]
which we could also write as:
\[\xy
(0,0)*+{\mbox{$\big[\cat{$(n-i-1)$-Gph}\big]$-\cat{Cat}-$\big[\cat{$i$-Gph}\big]$}
}="1";
(47,0)*+{\cat{$n$-Gph}}="2";
(33,0)*{{\scriptstyle \perp}};
{\ar@<1ex>_{} "1";"2"};
{\ar@<1ex>^{} "2";"1"};
\endxy\]
We write the induced monad as $T^{(n)}_i$, and its action is given by
\[T_iA(m) = \left\{
\begin{array}{lc}
A(m) & m \leq i \\[2ex]
\displaystyle\coprod_{k \geq 0}\quad \underbrace{A(m) \times_{A(i)} \cdots \times_{A(i)} A(m)}_{k \mbox{ {\small\it times}}} \hspace{2ex}& m>i 
\end{array}\right .\]
\end{prop}

This formula produces $k$-length strings of $i$-composable $m$-cells; as before we are taking $k$-fold wide pullbacks
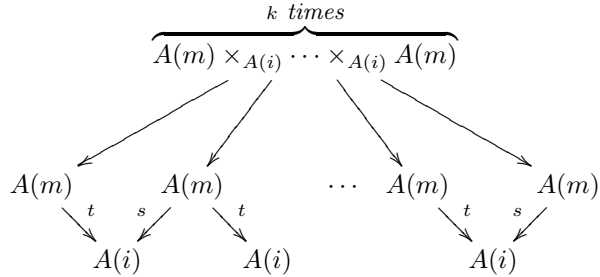
\[\xy
(0,10)*+{A(m)}="1";
(10,0)*+{A(i)}="2";
(20,10)*+{A(m)}="3";
(30,0)*+{A(i)}="4";
(50,10)*+{A(m)}="5";
(60,0)*+{A(i)}="6";
(70,10)*+{A(m)}="7";
(35,30)*+{\overbrace{A(m) \times_{A(i)} \cdots \times_{A(i)} A(m)}^{k \mbox{ {\small\it times}}}
}="8";
(40,10)*+{\cdots};
{\ar^{t} "1";"2"};
{\ar_{s} "3";"2"};
{\ar^{t} "3";"4"};
{\ar^{t} "5";"6"};
{\ar_{s} "7";"6"};
{\ar^{} "8";"1"};
{\ar_{} "8";"3"};
{\ar^{} "8";"5"};
{\ar^{} "8";"7"};
\endxy\]
where now $s$ and $t$ denote the composites along the top and bottom of 
\[\xy
(80,0)*+{A(i)}="2"; 
(63,0)*+{}="3"; 
(50,0)*+{}="4"; 
(30,0)*+{A(m-1)}="5"; 
(5,0)*+{A(m)}="6";
(57,0)*{\cdots};
%
%
%
{\ar@<0.7ex>^<<<<<<{s} "3"; "2"};
{\ar@<-0.7ex>_<<<<<<{t} "3"; "2"};
{\ar@<0.7ex>^>>>>>>>{s} "5"; "4"};
{\ar@<-0.7ex>_>>>>>>>{t} "5"; "4"};
{\ar@<0.7ex>^<<<<<<<{s} "6"; "5"};
{\ar@<-0.7ex>_<<<<<<<{t} "6"; "5"};
\endxy\]

\begin{proof}
By induction over $n$ and $i$.  Put $T^{(n)}_0 = \cat{fc}_{\cat{{\small ${\scriptstyle (n-1)}$-Gph}}}$, and for $i>0$ put $T^{(n)}_{i-1} = (T^{(n-1)}_i)_*$.
\end{proof}

We now show how to construct the distributive laws we require.  We will use the following proposition of Leinster; in fact this is just part of Proposition~F.1.1 of \cite{lei8}.  The notation may seem austere, but we will immediately give a motivating example below. 

\begin{prop}[Leinster \cite{lei8}]\label{key}
Let \cl{V} be a presheaf category and $T$ a monad on \cl{V}.  Write $\cl{V}^T$ for the category of algebras of $T$.  Then we have monads on \cat{$\cl{V}$-Gph} given by $T_*$ and $\cat{fc}_\cl{V}$, and a distributive law
\[ \lambda: T_* \circ \cat{fc}_\cl{V} \Rightarrow \cat{fc}_\cl{V} \circ T_*\]
whose resulting composite monad $\cat{fc}_\cl{V} \circ T_*$ is the free \cat{$\cl{V}^T$-Cat} monad, that is
\[ (\cat{\cl{V}-Gph})^{\cat{{\small fc}}_{\cl{V}} \circ T_*} \cong \cat{$\cl{V}^T$-Cat}.\]

\end{prop}

\begin{proof}
We have
\[ (T_* \circ \cat{fc}_\cl{V})  (A) = \coprod_{\begin{array}{c} {\scriptstyle k\geq 0} \\ {\scriptstyle a_0, \ldots, a_k \in A(0)} \end{array}} T\big(\ A(a_{k-1}, a_k) \times \cdots \times A(a_0, a_1)\ \big)\]
and
\[ (\cat{fc}_\cl{V} \circ T_*) (A) = \coprod_{\begin{array}{c} {\scriptstyle k\geq 0} \\ {\scriptstyle a_0, \ldots, a_k \in A(0)} \end{array}} T\big(A(a_{k-1}, a_k)\big) \times \cdots \times T\big(A(a_0, a_1)\big).\]
Now the universal property of the product 
\[T\big(A(a_{k-1}, a_k)\big) \times \cdots \times T\big(A(a_0, a_1)\big)\]
induces a canonical morphism from 
\[T\big(\ A(a_{k-1}, a_k) \times \cdots \times A(a_0, a_1)\ \big)\]
and this gives us the components of a natural transformation $\lambda$ as required.  It is straightforward to check that $\lambda$ is a distributive law.  
\end{proof}

The following example is a ``prototype'' for the construction of the strict $n$-category monad for general $n$.

\begin{eg}\label{revisit} {\bfseries 2-categories revisited.}

\noindent This time we put $\cl{V} = \cat{Gph}$ and $T = $ free category monad. Then we have

\begin{itemize}

\item $\cl{V}^T = \cat{Cat}$, 

\item $\cat{\cl{V}-Gph} = \cat{2-Gph}$,

\item $T_*$ is the monad on \cat{2-Gph} induced by the adjunction described in Example~\ref{ex1}, forming 1-composites of 2-cells freely, 

\item $\cat{fc}_{\cl{V}}$ is the monad on \cat{2-Gph} making free 0-composites, and 

\item $\lambda$ is given by the usual middle 4 interchange law for 2-categories.

\end{itemize}

\noindent The composite monad $\cat{fc}_\cl{V} \circ T_*$ resulting from this distributive law is the free 2-category monad on \cat{2-Gph}.  By the theory of distributive laws (Theorem~\ref{barrwells}) we also get a lift of the monad $\cat{fc}_{\cl{V}}$ to $T_*\mbox{-Alg} = \cat{Cat-Gph}$, whose algebras are precisely 2-categories.  This expresses 2-categories as graphs enriched in categories, with certain extra composition structure, which in effect gives us the usual definition of a 2-category as a category enriched in categories.

\end{eg}

\begin{eg}\label{ncatlei}{\bfseries $n$-categories as constructed by Leinster.}

\noindent In this example we recall Leinster's construction of the monad for strict $n$-categories, which is given as part of Theorem F.2.1 of \cite{lei8}.  The construction proceeds by induction.  We construct for each $n \geq 1$ a monad $S_n$ on \cat{$n$-Gph}, whose algebras are precisely strict $n$-categories.  We begin by taking $S_1$ to be the usual free category monad on \cat{Gph}.  Then for all $n \geq 2$ we apply Proposition~\ref{key} with

\begin{itemize}

\item $\cl{V} = \cat{$(n-1)$-Gph}$, and 

\item $T = S_{n-1}$, the free $(n-1)$-category monad that we have constructed by induction.

\end{itemize}

\noindent Then we have

\begin{itemize}

\item $\cl{V}^T = \cat{$(n-1)$-Cat}$, and

\item $\cat{\cl{V}-Gph} = \cat{$n$-Gph}$,

\end{itemize}
and the composite monad $\cat{fc}_{\cl{V}} \circ T_*$ resulting from the distributive law given by the Proposition has as its category of algebras

\[\cat{$\cl{V}^T$-Cat} = \cat{($(n-1)$-Cat)-Cat} = \cat{$n$-Cat}\]
which is to say that we have indeed constructed the free strict $n$-category monad.  As in the 2-category example above, we have essentially expressed $n$-categories as graphs enriched in $(n-1)$-categories, together with certain extra composition structure, which in effect gives us the usual definition of an $n$-category as a category enriched in $(n-1)$-categories.

\end{eg}

We now have everything we need to form all the distributive laws for interchange and thereby construct the monad for strict $n$-categories --- we simply start with a special case of Proposition~\ref{key} and then apply the $(-)_*$ construction repeatedly.  This is the content of Theorem~\ref{fourpointseven} and its proof.

\begin{thm}\label{fourpointseven}
The monads $T^{(n)}_0, \cdots, T^{(n)}_{n-1}$ on \cat{$n$-GSet} form a distributive series of monads as in Theorem~\ref{mainthm}.  For all $n>i>j\geq 0$ the distributive law
\[\lambda^{(n)}_{ij}: T^{(n)}_i T^{(n)}_j \Rightarrow T^{(n)}_j T^{(n)}_i\]
is given by interchange.  The resulting composite monad $T^{(n)}_0 T^{(n)}_1 \cdots T^{(n)}_{n-1}$ is the free strict $n$-category monad on \cat{$n$-GSet}.  

\end{thm}

\begin{proof}

First we construct the distributive laws.  Note that this, and indeed this whole proof, can be done directly by writing down and examining the formulae.  However we will take the more abstract approach. 

We begin by examining the case $j=0$, so we seek a distributive law
\[ T^{(n)}_i \circ T^{(n)}_0 \Rightarrow T^{(n)}_0 \circ T^{(n)}_i \]
for each $n>i>0$.  But we know 
\[T^{(n)}_i = (T^{(n-1)}_{i-1})_*\] 
and 
\[T^{(n)}_0 = \cat{fc}_{\cat{{\small ${\scriptstyle (n-1)}$-Gph}}}\]
so this is just a special case of Proposition~\ref{key} above.  

For $j>0$ we use the distributive law 
\[ T^{(n-j)}_{i-j} \circ T^{(n-j)}_0 \Rightarrow T^{(n-j)}_0 \circ T^{(n-j)}_{i-j} \]
and apply the $(-)_*$ construction $j$ times; since this is a 2-functor, the result of applying it to a distributive law must be a distributive law.  

Finally it is straightforward to check that these distributive laws obey all the necessary Yang-Baxter equations.  

To show that the composite monad $T^{(n)}_0 T^{(n)}_1 \cdots T^{(n)}_{n-2} T^{(n)}_{n-1}$ is the free strict $n$-category monad we also proceed by induction.  The result is clearly true for $n=1$.  For $n>1$ we know by Theorem~\ref{mainthm} that the this composite monad arises from various different composite distributive laws; in particular it arises from the distributive law

\[ \big( T^{(n)}_1 T^{(n)}_2 \cdots T^{(n)}_{n-1}\big) \circ T^{(n)}_0 \Rightarrow T^{(n)}_0 \circ \big( T^{(n)}_1 T^{(n)}_2  \cdots T^{(n)}_{n-1} \big) .\]
Now by definition we have
\begin{eqnarray*}
T^{(n)}_1 T^{(n)}_2 \cdots T^{(n)}_{n-1} &=& (T^{(n-1)}_0)_* (T^{(n-1)}_1)_* \cdots (T^{(n-1)}_{n-2})_* \\
&=& \big( T^{(n-1)}_0 T^{(n-1)}_1 \cdots T^{(n-1)}_{n-2} \big)_* 
\end{eqnarray*}

\noindent but by induction 
\[T^{(n-1)}_0 T^{(n-1)}_1 \cdots T^{(n-1)}_{n-2}\]
is the free strict $(n-1)$-category monad.  So this distributive law is exactly the one that Leinster uses to construct the free strict $n$-category monad.  
\end{proof}

We will now illustrate this construction for the case $n=2, j=0, i=1$, which should be the usual interchange law between horizontal and vertical composition.

\begin{itemize}

\item A cell of $T^{(2)}_0(A)$ is an $l$-length string of 0-composable cells, for example
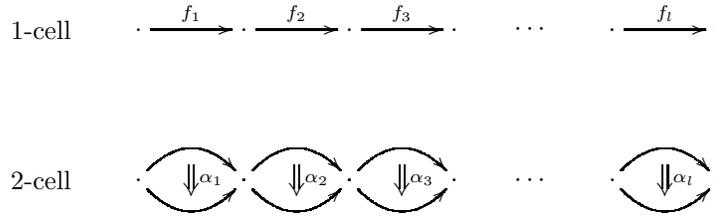
\[\xy
(-20,20)*{\mbox{1-cell}};
(-20,0)*{\mbox{2-cell}};
(-7,20)*+{.}="1";
(7,20)*+{.}="2";
(21,20)*+{.}="3";
(35,20)*+{.}="4";
(56,20)*+{.}="5";
(70,20)*+{.}="6";
(45,20)*{\cdots};
{\ar^{f_1} "1";"2"};
{\ar^{f_2} "2";"3"};
{\ar^{f_3} "3";"4"};
{\ar^{f_l} "5";"6"};
(0,0)*{\xy
(-7,0)*+{.}="1";
(7,0)*+{.}="2";
{\ar@/^1pc/^{} "1";"2"};
{\ar@/_1pc/_{} "1";"2"};
{\ar@{=>}^{{\scriptstyle \alpha_1}} (0,2)*{};(0,-2)*{}} ;
\endxy};
(14,0)*{\xy
(-7,0)*+{.}="1";
(7,0)*+{.}="2";
{\ar@/^1pc/^{} "1";"2"};
{\ar@/_1pc/_{} "1";"2"};
{\ar@{=>}^{\alpha_2} (0,2)*{};(0,-2)*{}} ;
\endxy};
(28,0)*{\xy
(-7,0)*+{.}="1";
(7,0)*+{.}="2";
{\ar@/^1pc/^{} "1";"2"};
{\ar@/_1pc/_{} "1";"2"};
{\ar@{=>}^{\alpha_3} (0,2)*{};(0,-2)*{}} ;
\endxy};
(63,0)*{\xy
(-7,0)*+{.}="1";
(7,0)*+{.}="2";
{\ar@/^1pc/^{} "1";"2"};
{\ar@/_1pc/_{} "1";"2"};
{\ar@{=>}^{\alpha_l} (0,2)*{};(0,-2)*{}} ;
\endxy};
(45,0)*{\cdots};
\endxy\]

\item A cell of $T^{(2)}_1(A)$ is an $h$-height ``column'' of 1-composable cells, for example a 2-cell has the form
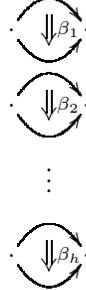
\[\xy
(0,0)*{\xy
(-5,0)*+{.}="1";
(5,0)*+{.}="2";
{\ar@/^1pc/^{} "1";"2"};
{\ar@/_1pc/_{} "1";"2"};
{\ar@{=>}^{{\scriptstyle \beta_h}} (0,2)*{};(0,-2)*{}} ;
\endxy};
(0,20)*{\xy
(-5,0)*+{.}="1";
(5,0)*+{.}="2";
{\ar@/^1pc/^{} "1";"2"};
{\ar@/_1pc/_{} "1";"2"};
{\ar@{=>}^{\beta_2} (0,2)*{};(0,-2)*{}} ;
\endxy};
(0,30)*{\xy
(-5,0)*+{.}="1";
(5,0)*+{.}="2";
{\ar@/^1pc/^{} "1";"2"};
{\ar@/_1pc/_{} "1";"2"};
{\ar@{=>}^{\beta_1} (0,2)*{};(0,-2)*{}} ;
\endxy};
(0,11)*{\vdots};
\endxy\]

\item A 2-cell of $T^{(2)}_1 T^{(2)}_0(A)$ is thus an 
\[\mbox{``$h$-height column of 1-composable \{~strings~of~0-composable~2-cells~\}''.}\]

  Note that for these strings of 0-composable cells to be 1-composable, they must all have the same length $l$.  So we have an $l \times h$ composable ``grid'' of 2-cells
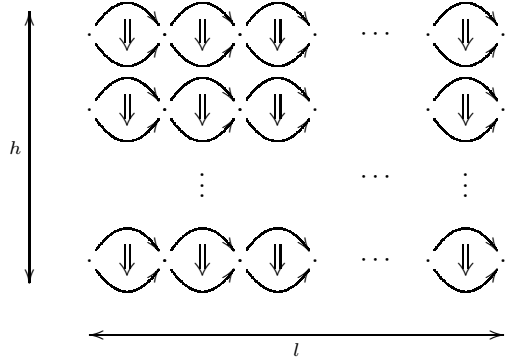
\[\xy
(0,0)*{\glob};
(10,0)*{\glob};
(20,0)*{\glob};
(33,0)*{\cdots};
(45,0)*{\glob};
(0,20)*{\glob};
(10,20)*{\glob};
(20,20)*{\glob};
(45,20)*{\glob};
(0,30)*{\glob};
(10,30)*{\glob};
(20,30)*{\glob};
(33,30)*{\cdots};
(45,30)*{\glob};
(10,11)*{\vdots};
(45,11)*{\vdots};
(33,11)*{\cdots};
{\ar@{<->}_{l} (-5,-10)*{};(50,-10)*{}};
{\ar@{<->}^{h} (-13,-3)*{};(-13,33)*{}};
\endxy\]

\item On the other hand a 2-cell of $T^{(2)}_0 T^{(2)}_1(A)$ is an 
\[\mbox{``$l$-length string of 0-composable \{columns of 1-composable 2-cells}''.\]

\noindent  Note that for the columns to be 0-composable they do not have to be the same height, so we have a configuration of 2-cells such as:
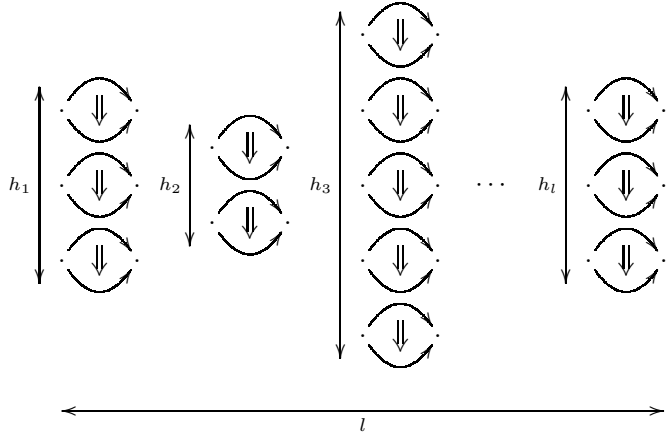
\[\xy
(0,-10)*{\glob};
(0,0)*{\glob};
(0,10)*{\glob};
{\ar@{<->}^{h_1} (-8,-13)*{};(-8,13)*{}};
(20,-5)*{\glob};
(20,5)*{\glob};
{\ar@{<->}^{h_2} (12,-8)*{};(12,8)*{}};
(40,-20)*{\glob};
(40,-10)*{\glob};
(40,0)*{\glob};
(40,10)*{\glob};
(40,20)*{\glob};
{\ar@{<->}^{h_3} (32,-23)*{};(32,23)*{}};
(70,-10)*{\glob};
(70,0)*{\glob};
(70,10)*{\glob};
{\ar@{<->}^{h_l} (62,-13)*{};(62,13)*{}};
{\ar@{<->}_{l} (-5,-30)*{};(75,-30)*{}};
(52,0)*{\cdots};
\endxy\]

\end{itemize}

There is evidently a natural map from an $l \times h$ grid to one with columns of height $h_1, \cdots, h_l$ as above -- we simply put $h_1 = h_2 = \cdots = h_l = h$.  This is the canonical map
\[T\big(\ A(a_{l-1}, a_l) \times \cdots \times A(a_0, a_1)\ \big) \lra
T\big(A(a_{l-1}, a_l)\big) \times \cdots \times T\big(A(a_0, a_1)\big)\]
as described in the proof of Proposition~\ref{key}.  In our case we are using $T_* = T^{(2)}_1$ so $T = T^{(1)}_0$ which is just the free category monad.  Note that the left hand side \[T\big(\ A(a_{l-1}, a_l) \times \cdots \times A(a_0, a_1)\ \big)\] is $T$ applied to the product \[A(a_{l-1}, a_l) \times \cdots \times A(a_0, a_1).\]  We can express this product as a (trivial) wide pullback
%
%
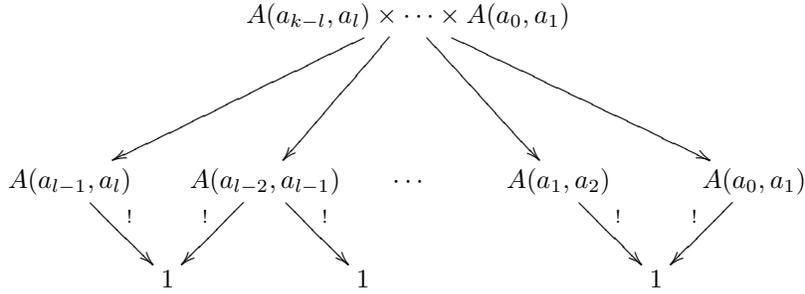
\[\xy
(0,13)*+{A(a_{l-1}, a_{l})}="1";
(13,0)*+{1}="2";
(26,13)*+{A(a_{l-2}, a_{l-1})}="3";
(39,0)*+{1}="4";
(65,13)*+{A(a_1, a_2)}="5";
(78,0)*+{1}="6";
(91,13)*+{A(a_0, a_1)}="7";
(45,35)*+{A(a_{k-l}, a_l) \times \cdots \times A(a_0, a_1)}="8";
(45,13)*+{\cdots};
{\ar^{!} "1";"2"};
{\ar_{!} "3";"2"};
{\ar^{!} "3";"4"};
{\ar^{!} "5";"6"};
{\ar_{!} "7";"6"};
{\ar^{} "8";"1"};
{\ar_{} "8";"3"};
{\ar^{} "8";"5"};
{\ar^{} "8";"7"};
\endxy\]
and since $T$ is cartesian we know that applying it to this product gives the wide pullback
%
%
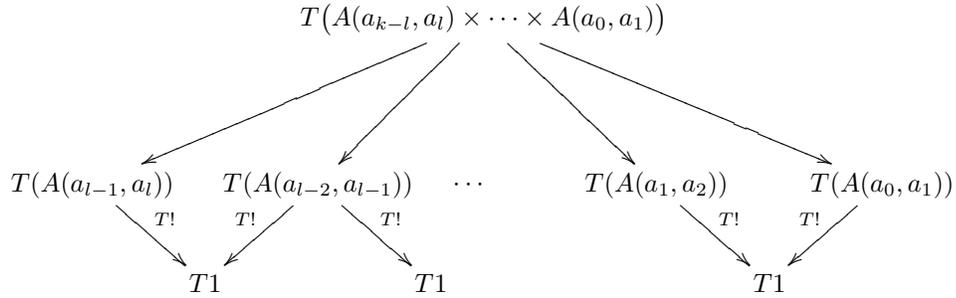
\[\xy
(0,13)*+{T(A(a_{l-1}, a_{l}))}="1";
(15,0)*+{T1}="2";
(30,13)*+{T(A(a_{l-2}, a_{l-1}))}="3";
(45,0)*+{T1}="4";
(75,13)*+{T(A(a_1, a_2))}="5";
(90,0)*+{T1}="6";
(105,13)*+{T(A(a_0, a_1))}="7";
(52,35)*+{T\big(A(a_{k-l}, a_l) \times \cdots \times A(a_0, a_1)\big)}="8";
(50,13)*+{\cdots};
{\ar^{T!} "1";"2"};
{\ar_{T!} "3";"2"};
{\ar^{T!} "3";"4"};
{\ar^{T!} "5";"6"};
{\ar_{T!} "7";"6"};
{\ar^{} "8";"1"};
{\ar_{} "8";"3"};
{\ar^{} "8";"5"};
{\ar^{} "8";"7"};
\endxy\]
The commuting condition over $T1$ is what ensures that all the ``columns'' must now have equal height.  

On the other hand there is not a distributive law going in the opposite direction, since given a grid with columns of possibly varying height, there is no canonical way to map it to a grid with columns of equal height.  We might attempt to insert 2-cell identities to ``extend'' the shorter columns, but the possible choice of positions for the inserted identity cells means that this will not satisfy the axioms for a distributive law.


\begin{thebibliography}{10}

\bibitem{bw1}
M.~Barr and C.~Wells.
\newblock {\em Toposes, Triples and Theories}.
\newblock Number 278 in Grundlehren der math. Wissenschaften. Springer-Verlag,
  1985.

\bibitem{bat1}
M.~A. Batanin.
\newblock Monoidal globular categories as a natural environment for the theory
  of weak $n$-categories.
\newblock {\em Adv. Math.}, 136(1):39--103, 1998.

\bibitem{bec1}
J.~Beck.
\newblock Distributive laws.
\newblock {\em Lecture Notes in Mathematics}, 80:119--140, 1969.

\bibitem{che16}
Eugenia Cheng.
\newblock Operadic theories of $n$-category, 2007.
\newblock Talk at CT07, Carvoeiro.

\bibitem{gps1}
R.~Gordon, A.~J. Power, and R.~Street.
\newblock Coherence for tricategories.
\newblock {\em Memoirs of the American Mathematical Society}, 117(558), 1995.

\bibitem{gra1}
J.~W. Gray.
\newblock {\em Formal category theory: adjointness for 2-categories}, volume
  391 of {\em Lecture Notes in Mathematics}.
\newblock Springer-Verlag, Berlin-New York, 1974.

\bibitem{gra2}
J.~W. Gray.
\newblock Coherence for the tensor product of 2-categories, and braid groups.
\newblock In {\em Algebra, Topology and Category Theory, a collection in honor
  of Samuel Eilenberg}, LMS, pages 63--76. Academic Press, 1976.

\bibitem{gur2}
Nick Gurski.
\newblock {\em An algebraic theory of tricategories}.
\newblock PhD thesis, University of Chicago, June 2006.
\newblock Available via {\tt http://www.math.yale.edu/$\sim$mg622/tricats.pdf}.

\bibitem{jk1}
Andr\'e Joyal and Joachim Kock.
\newblock Weak units and homotopy 3-types, 2005.
\newblock To appear in \emph{Categories in Algebra, Geometry and Mathematical
  Physics}, proceedings of Streetfest 2005, eds Batanin, Davydov, Johnson,
  Lack, Neeman; preprint available via {\tt math.CT/0602084}.

\bibitem{lei8}
Tom Leinster.
\newblock {\em Higher operads, higher categories}.
\newblock Number 298 in London Mathematical Society Lecture Note Series.
  Cambridge University Press, 2004.

\bibitem{sim2}
Carlos Simpson.
\newblock A closed model structure for $n$-categories, internal {H}om,
  $n$-stacks and generalized {S}eifert-{V}an~{K}ampen, 1997.
\newblock Preprint alg-geom/9704006.

\bibitem{sim1}
Carlos Simpson.
\newblock Limits in $n$-categories, 1997.
\newblock Preprint {\tt alg-geom/9708010}.

\bibitem{str1}
Ross Street.
\newblock The formal theory of monads.
\newblock {\em Journal of Pure and Applied Algebra}, 2:149--168, 1972.

\bibitem{str2}
Ross Street.
\newblock The algebra of oriented simplexes.
\newblock {\em Journal of Pure and Applied Algebra}, 49(3):283--335, 1987.

\bibitem{tri1}
Todd Trimble.
\newblock What are `fundamental $n$-groupoids'?, 1999.
\newblock seminar at DPMMS, Cambridge, 24 August 1999.

\end{thebibliography}

\ed